\documentclass[11pt]{article}

\usepackage[margin=1in]{geometry}

\usepackage{times}
\usepackage{url}
\usepackage{fullpage}
\usepackage{amssymb}
\usepackage{xspace}
\usepackage{amsthm}
\usepackage{cite}
\usepackage{caption}
\usepackage{amsmath}
\usepackage{amsfonts}
\usepackage{times}
\usepackage{epstopdf}
\usepackage{amsmath}
\usepackage{color}
\usepackage{amsmath}
\usepackage{float}
\usepackage{defs}
\usepackage{clrscode}
\usepackage{textlater}
\usepackage{paralist}
\usepackage{tocbibind}
\usepackage{graphicx}
\usepackage{epstopdf}
\usepackage{epsfig}
\usepackage{xspace}
\usepackage{amsmath}
\usepackage{amssymb}
\usepackage{amsthm}
\usepackage{amsfonts}
\usepackage{array}
\usepackage{url}
\usepackage{multirow}
\usepackage{pifont}
\usepackage{bbm}
\usepackage{lineno}
\usepackage{times}
\usepackage{color}
\usepackage{wrapfig}
\usepackage{float}
\usepackage{cite}
\usepackage{epsfig}
\usepackage{hyphenat}
\usepackage[stable]{footmisc}
\usepackage[T1]{fontenc}
\usepackage{dcolumn}
\usepackage{endnotes}
\usepackage{graphics}
\usepackage{graphicx}
\usepackage{framed}
\usepackage{verbatim}
\usepackage{placeins}
\usepackage{longtable}
\newtheorem{problem}[theorem]{Problem}

\newtheorem{program}[theorem]{Program}

\def\R{{\cal R}}

\newcommand{\al}{\alpha}
\newcommand{\be}{\beta}
\newcommand{\ga}{\gamma}

\newenvironment{CompactItemize}
 { \begin{itemize}}
 {\end{itemize}}
\newenvironment{CompactEnumerate}
 { \begin{enumerate}}
 {\end{enumerate}}

\newenvironment{proofof}[1]{\noindent{\bf Proof of {#1}.}  \ignorespaces}{\null\hfill\qedsymbol}

\usepackage{full}

\begin{document}

\begin{titlepage}
%\title{Extremal Quasiconformal Maps between Planar Domains}
%
%
%\author{
%  Mayank Goswami\thanks{ Department of Applied Mathematics and Statistics 
%Stony Brook University, 
%Stony Brook, NY 11794-3600 USA.
%Email: \texttt{mgoswami@ams.sunysb.edu}.}
%\and
%Xianfeng Gu\thanks{Dept.\ of Computer Science, Stony Brook University, Stony Brook, NY 11791-4400 USA.
%Email: \texttt{\{gu}@cs.stonybrook.edu}.}
%}
%
%\iffalse
%
%\thanks{Stony Brook University\thanks{\sloppy \email{\{gu}@cs.stonybrook.edu},
%  \email{mgoswami@ams.sunysb.edu}.}
%}

%!TEX root =  ak-sorting.tex

%%%%%%%%%%%%%%
\iffalse
%%%%%%%%%%%%%%

MISTRESS QUICKLY: Marry, this is the short and the long of it; you have brought her into such a canaries as 'tis wonderful. The best courtier of them all, when the court lay at Windsor, could never have brought her to such a canary

The Long and the S~ort of It: External Memory Sorting with Different Length Atomic Keys
\title{Sorting Shoes, Ships, Cabbages, and Kings: External Memory Sorting with Different Length Atomic Keys}

%%%%%%%%%%%%%%
\fi
%%%%%%%%%%%%%%

\title{Computing Teichm\"{u}ller Maps between Polygons}

\author{
  Mayank Goswami\thanks{Max Planck Institute for Informatics,
  Saarbr\"{u}cken 66123, Germany.
Email: \texttt{gmayank@mpi-inf.mpg.de}. Part of this work was done when the author was a graduate student in the Department of Applied Mathematics and Statistics at Stony Brook University.}
% MAB: actually I don't know which dept Mayank is in. 
  \and 
 Xianfeng Gu\thanks{Department of Computer Science, 
 Stony Brook University,
 Stony Brook, NY 11794-4400, USA.
 Email: \texttt{gu@cs.stonybrook.eu}.}
  \and 
 Vamsi P. Pingali \thanks{Department of Mathematics, 
Johns Hopkins University,
Baltimore, MD - 21218, USA.
 Email: \texttt{vpingali@math.jhu.edu}.}
 \and
 Gaurish Telang \thanks{Department of Applied Mathematics and Statistics,
 Stony Brook University,
 Stony Brook, NY 11794, USA.
 Email: \texttt{gaurish.ams@gmail.com}
 }
}

\iffalse

\thanks{Stony Brook University\thanks{\sloppy \email{gu\@cs.stonybrook.edu},
  \email{mgoswami@ams.sunysb.edu}.}

}
\fi

\date{}

\maketitle \thispagestyle{empty}

\begin{abstract}
%
%By the Riemann-mapping theorem, one can always bijectively map an $n$-polygon $P$ to another $n$ polygon $Q$ conformally; in fact, any two simply connected domains can be conformally mapped to each other. However, this mapping need not necessarily map the vertices of $P$ to those $Q$. In this case, one wants to find the ``best" mapping between these polygons, i.e., one that minimizes the maximum angle distortion over all points in the base polygon. This polygon mapping problem is a special instance of the extremal quasiconformal map problem on punctured spheres.
%
%While existence and characterization of extremal quasiconformal maps is well studied in mathematics, in this paper we provide the first explicit construction of the unique extremal map in the continuous setting. Our construction is via an iterative procedure which is proven to converge to the unique extremal map.
%
%We also provide the first numerical method to solve the extremal quasiconformal map problem on punctured spheres. As a result, the algorithm for the polygon mapping problem computes the image of a dense sample of points from the polygon under the extremal quasiconformal map. Our method uses a variational approach to successively improve the maximum angle distortion of an initial map.

By the Riemann-mapping theorem, one can bijectively map the interior of an $n$-gon $P$ to that of another $n$-gon $Q$ conformally. However, (the boundary extension of) this mapping need not necessarily map the vertices of $P$ to those $Q$. In this case, one wants to find the ``best" mapping between these polygons, i.e., one that minimizes the maximum angle distortion (the dilatation) over \textit{all} points in $P$. From complex analysis such maps are known to exist and are unique. They are called extremal quasiconformal maps, or Teichm\"{u}ller maps.

Although there are many efficient ways to compute or approximate conformal maps, there is currently no such algorithm for extremal quasiconformal maps. This paper studies the problem of computing extremal quasiconformal maps both in the continuous and discrete settings.

We provide the first constructive method to obtain the extremal quasiconformal map in the continuous setting. Our construction is via an iterative procedure that is proven to converge quickly to the unique extremal map. To get to within $\epsilon$ of the dilatation of the extremal map, our method uses $O(1/\epsilon^{4})$ iterations. Every step of the iteration involves convex optimization and solving differential equations, and guarantees a decrease in the dilatation. Our method uses a reduction of the polygon mapping problem to that of the punctured sphere problem, thus solving a more general problem.

We also discretize our procedure. We provide evidence for the fact that the discrete procedure closely follows the continuous construction and is therefore expected to converge quickly to a good approximation of the extremal quasiconformal map.

\end{abstract} 

\end{titlepage}
%\linenumbers

\section{Introduction}
\label{sec:introduction}

One of the foundational results in complex analysis, the Riemann mapping theorem, states that any non-empty simply connected domain $U \subsetneq \mathbb{C}$ can be mapped bijectively and conformally to the unit disk $\mathbb{D}$. This implies that the interiors of two simple planar $n$-gons $P$ and $Q$ can be mapped bijectively and conformally to each other. By another result \cite{cara}, such a map $f: P \rightarrow Q$ extends continuously to the boundary $\bar{P}$ of $P$ (the edges). Generally, the vertices of $P$ do not map to the vertices of $Q$ under this extended mapping.

Assume we are given an ordering $\{v_{i}\}_{i=1}^{n}$ and $\{v_{i}^{'}\}_{i=1}^{n}$ of the vertices of $P$ and $Q$, respectively. Consider the space of homeomorphisms  $f$ that map $\bar{P}$ to $\bar{Q}$, such that $f(v_{i})=v_{i}^{'}$. Such an $f$ is bound to stretch angles (unless the polygons are linear images of each other), and a classical way to measure this angle stretch by $f$ at a point $p \in P$ is by $\mu_{f}(p) = f_{\bar{z}}(p) / f_{z}(p)$. This complex-valued function $\mu_{f}$ is called the Beltrami coefficient of $f$, and it satisfies $||\mu_{f}||_{\infty} <1$. The problem we consider is computing the "best" homeomorphism $f_*$ in the above class, i.e., an $f_*$ such that the norm of its Beltrami differential $||\mu_*||_{\infty}$ is the smallest amongst all homeomorphisms satisying the above conditions. These homeomorphisms that stretch angles but by a bounded amount are called quasiconformal homeomorphisms (q.c.h.), and the "best" q.c.h. $f_*$ is called the extremal quasiconformal map, or the Teichm\"{u}ller map.

As an example, consider two rectangles $R_{i} = [0,a_{i}] \times [0,b_{i}](i = 1,2)$ in the plane. Starting from the origin, label the vertices of $R_{1}$ and $R_{2}$ counter-clockwise as $\{v_{j}\}_{j=1}^{4}$ and $\{v_{j}^{'}\}_{j=1}^{4}$, respectively ($v_{1}=v_{1}^{'}=(0,0)$). Consider the space of all q.c.h. $f:R_{1} \rightarrow R_{2}$ such that $f(v_{i}) = v_{i}^{'}$. It was shown by Gr\"{o}tzsch \cite{gro} that the affine map $f_*(x,y) = (a_{2}x/a_{1}, b_{2}y/b_{1})$ is the unique extremal quasiconformal map; any other map $f$ would stretch angles at some point $p \in R_{1}$ more than $g$ (i.e., $\exists p \in R_{1}:|\mu_{f}(p)| > |\mu_{*}(p)|$). For the general $n$-gon case mentioned above, such an analytic solution does not exist. However, the extremal map exists and is unique (these are the famous theorems by Teichmuller \cite{teich1} and \cite{teich2}, proven rigorously later by Ahlfors \cite{Ahlfors_2006}), and is of the form $\text{conformal} \circ \text{affine} \circ \text{conformal}$.

Computing a Riemann mapping from a given polygon to the disk has gathered a lot of attention in the past. Algorithms (e.g., CRDT\cite{Driscoll_numericalconformal}) based on finding the unknown parameters in the Schwarz-Christoffel mapping formula \cite{driscoll2002schwarz} for a Riemann map were proposed, and the latest result by Bishop \cite{Bishoplinear} computes a $(1+\epsilon)$ quasiconformal map in $O(n\log (1/\epsilon) \log \log (1/\epsilon))$.

No such algorithm that computes (or approximates) the extremal quasiconformal map is known. In contrast to the Riemann mapping problem, where Riemann gave a constructive proof, the proof by Teichm\"{u}ller/Ahlfors is an existence result, and no constructive proofs are available. Furthermore, the "formula" for Teichm\"{u}ller maps analogous to the Schwarz-Christoffel mapping for Riemann mapping states that $\mu_{*} = k \bar{\phi} / |\phi|$, for some integrable holomorphic function $\phi$, with at most simple poles at the vertices of $P$. Thus, given $P$, we know all the extremal maps with domain $P$; our problem is figuring out which one takes us to our target $Q$. Even though $\phi$ comes from a finite dimensional family, there is no direct search criterion\footnote{It is not known how much a variation is $\phi$ would change the solution of the Beltrami equation for $k \bar{\phi} / |\phi|$}. This should be contrasted with the known relation between the images of the vertices of the polygon in the Schwarz-Christoffel formula and the concept of harmonic measure \cite{HarmonicMeasure05Book}. In fact, to the authors' knowledge, there does not exist a method that, given a starting $f$ between $P$ and $Q$, computes a $g$ with $||\mu_{g}||_{\infty} < ||\mu_{f}||_{\infty}$ if one exists.

This paper gives the first results for theoretically constructing and algorithmically computing Teichm\"{u}ller maps for the polygon case stated above. Our procedure is iterative; we 1) start with a q.c.h. that sends the vertices of $P$ to the vertices of $Q$ in the prescribed order, 2) improve on it, and then 3) recurse on the improved map.

The problem of computing a Teichm\"{u}ller map is syonymous with
computing geodesics in the Teichm\"{u}ller space endowed with the
Teichm\"{u}ller metric, which is the universal cover of the moduli
space of Riemann surfaces (in which all (mutually homeomorphic)
one-dimensional complex manifolds are quotiented under the
equivalence relation of biholomorphism). Teichm\"{u}ller theory is
an active area of research in mathematics, and it has connections to
topology\footnote{It has been used by Lipman Bers to give a simpler
proof of Thurston's classification theorem for surface
homeomorphisms.}, dynamics, algebraic geometry, and number theory.
Being able to compute the distance between two given points in a
Teichm\"{u}ller space (two equivalence classes of marked Riemann
surfaces) would help us learn more about the geometry of this
interesting space. This work is therefore intended to be an
introduction to this rich subject from a computational perspective,
and we certainly feel that many computationally challenging open
problems lie hidden.

Computing Teichm\"{u}ller maps is also an important problem in the fields of medical imaging, computer graphics and vision. In medical imaging, conformal and quasiconformal mapping has been applied for brain cortical surface registration (\cite{Wang:TMI}, \cite{Gu04}). In computer vision, conformal geometry has been applied for shape analysis and dynamic surface registration and tracking (\cite{DBLP:journals/ijcv/WangGZWGSH08}),\cite{TPAMI10Ricci}), and in computer graphics, conformal geometry has been applied for surface parameterization (\cite{Gu03SGP}).

Surface registration refers to the process of finding an optimal one-to-one correspondence between surfaces that preserves the surface geometric structures and reduces the distortions as much as possible. Teichm\"{u}ller maps satisfy all these requirements. Thus being able to compute them would help one get a novel algorithm for surface registration. In \cite{Zorin12SGP} various advantages of extremal quasiconformal maps over many existing methods were discussed in detail, and we refer the reader to it for an overview of how extremal quasiconformal maps are important in geometry processing.

\paragraph{Related work}
The only previous work to have considered the problem of computing extremal quasiconformal maps is \cite{Zorin12SGP}. The authors consider a very similar version where a Dirichlet boundary condition is given on the disk, and one is required to compute the extremal map whose boundary values satisfy the given condition. The authors propose a heuristic; they obtain a "highly nonlinear" energy and minimize it using an alternate-descent method. There is no guarantee on how far the solution is from the true extremal map, as the solution obtained could be a local minima of the energy. Another possibly related work is \cite{Sharon:2006:AUC}, where the authors use the concept of conformal welding to get fingerprints for a simple closed curve.

Various eminent mathematicians (Teichm\"{u}ller, Ahlfors, Bers, Reich, Strebel, Krushkal, Hamilton, etc.) have contributed to Teichm\"{u}ller theory. We refer the reader to \cite{Gardiner} and \cite{hubbard2006} for some excellent introductions to Teichm\"{u}ller theory. Most of the classical results we use can either be found in these books, or references contained therein.

\section{Problem statements and results}
\seclabel{pdef}

In this section we first state rigorously what extremal quasiconformal map we want to compute, and what we mean by computing such a map. We will then state our main results.

\subsection{Problem statements}

The amount of angle stretch induced by a quasiconformal homeomorphism (abbreviated henceforth as q.c.h.) $f$ can be quantified using the Beltrami coefficient $\mu_{f}$ of $f$. Defining $f_{\bar{z}} = f_{x} + i f_{y}$ and $f_{z}= f_{x} - i f_{y}$, where $f_{x}$ and $f_{y}$ denote the partials of $f$ w.r.t. $x$ and $y$, the Beltrami coefficient $\mu_{f}$ is defined as $\mu_{f}=f_{\bar{z}}/f_{z}$. Intuitively, a q.c.h. maps the unit circle in the tangent space at a point $p$ in the domain to an ellipse in the tangent space at $f(p)$, and $(1+|\mu_{f}|)/(1-|\mu_{f}|)$ is the eccentricity of this ellipse. For a formal definition of quasiconformal maps and Beltrami differentials, we refer the reader to \subsecref{qcmapsapp} in the Appendix.

Let $P$ and $Q$ be two $n$-polygons\footnote{We allow for $\infty$ to be a vertex of the polygon.} in the plane. Let $\{v_{i}\}_{i=1}^{n}$ and $\{v_{i}^{'}\}_{i=1}^{n}$ be an ordering of the vertices of $P$ and $Q$, respectively. Observe that: 
\begin{CompactEnumerate}
\item The polygons, or for that matter any simply connected domain (with boundary as a Jordan curve) is conformally equivalent to the upper half plane $\mathbb{H}$, and \item Composition by conformal maps does not change the dilatation (maximal angle stretch).
\end{CompactEnumerate}

Therefore, an $n$-gon is the same as $\mathbb{H}$ with $n$ marked points on the boundary $\partial \mathbb{H} = \mathbb{R}$.
%Since the polygons are conformally equivalent to the disk, each of the polygons can be conformally mapped to the unit disk $\mathbb{D}$. Assume that $z_{i}$ and $w_{i}$ are the images (under a conformal map) of the vertices of $P$ and $Q$, respectively. The above problem then translates to:

\begin{problem}\label{polygonmappingproblem}[Polygon mapping problem]
Given $\{z_{1},...z_{n}, w_{1},...w_{n}\} \in \partial \overline{\mathbb{H}}$, find $f_{*}:\overline{\mathbb{H}} \rightarrow \overline{\mathbb{H}}$ (with Beltrami coefficient $\mu_*$) satisfying:
\begin{CompactEnumerate}

\item $f_*$ is a quasiconformal homeomorphism of $\overline{\mathbb{H}}$ to itself.
\item $f_*(z_{i})=w_{i}$, $i \in \{1,...n\}$
\item $ || \mu_{*}||_{\infty} \leq ||\mu_{f}||_{\infty} $ for all $f$ satisfying (1) and (2) above.
\end{CompactEnumerate}
\end{problem}
Note that by Teichm\"{u}ller's theorems, the above $f_*$ exists and is unique. We state the punctured sphere problem next, and show that it is in fact a generalization of the polygon mapping problem.

\begin{problem}\label{puncturedsphereproblem}[Punctured sphere problem]
Given $\{z_{1},...z_{n-3},z_{n-2}=0,z_{n-1}=1,z_{n}=\infty\}$, $\{w_{1},...w_{n-3}$, $w_{n-2}=0,w_{n-1}=1,w_{n}=\infty\}$, and $h:\hat{\mathbb{C}} \rightarrow \hat{\mathbb{C}}$ such that $h(z_{i})=w_{i}$, find $f_*:\hat{\mathbb{C}} \rightarrow \hat{\mathbb{C}}$ satisfying:

\begin{CompactEnumerate}
\item $f_*$ is a quasiconformal homeomorphism of $\hat{\mathbb{C}}$ to itself.
\item  $f_*$ is isotopic to $h$ relative to the points $\{0,1,\infty,z_{1},..z_{n-3}\}$, i.e. $f_*(z_{i})=w_{i}$.
\item $ || \mu_{*}||_{\infty} \leq ||\mu_{f}||_{\infty} $ for all $f$ satisfying (1) and (2) above.
\end{CompactEnumerate}

\end{problem}

We call the base $z_{i}$-punctured sphere $R$ and the target $w_{i}$-punctured sphere $S$ from now on. The reason why the punctured sphere problem requires a starting map $h$ as input is that by Teichm\"{u}ller's theorem, the extremal map exists and is unique within each homotopy class. The following theorem shows that Problem~\ref{puncturedsphereproblem} is indeed general.

\begin{theorem}\thmlabel{problemsrelation}
An algorithm for Problem~\ref{puncturedsphereproblem} can be used to give a solution to Problem~\ref{polygonmappingproblem}.
\end{theorem}

\textlater{\begin{proofof}{\thmref{problemsrelation}}
We take an instance of the polygon mapping problem and conver it to an instance of the punctured sphere problem first.

Let $h_{0}$ be any quasiconformal homeomorphism mapping $P$ to $Q$, such that $h_{0}(v_{i}) = v_{i}^{'}$. By conformally mapping $P$ and $Q$ to $\mathbb{H}$ (denote the maps by $\pi_{P}$ and $\pi_{Q}$), we get a quasiconformal self-homeomorphism $h_{u}$ of $\mathbb{H}$, satisfying $h_{u}(z_{i})=w_{i}$, where $z_{i}$ and $w_{i}$ are images (under $\pi_{P}$ and $\pi_{Q}$) of $v_{i}$ and $v_{i}^{'}$, respectively. Furthermore $h_{u}$ can be normalized to fix $0$, $1$ and $\infty$. Let $\overline{\mathbb{H}}$ denote the lower half plane, and define a quasiconformal self-homeomorphism $h_{\ell}$ of $\overline{\mathbb{H}}$ by $h_{\ell}(z) = \overline{h_{u}(\bar{z})}$. Now $h_{u}$ and $h_{\ell}$ agree on $\mathbb{R}$, and can be pieced together to get a  quasiconformal self-homeomorphism $h$ of $\hat{\mathbb{C}}$ satisfying $h(z_{i})=w_{i}$. Note that $h$ fixes $0$, $1$ and $\infty$.

The next theorem shows how one can get back the answer to the polygon mapping problem from the answer to the punctured sphere problem.

\begin{lemma}\lemlabel{thmrelation}
Let $f$ be the solution to Problem~\ref{puncturedsphereproblem} when it is fed the input data $\{z_{i},w_{i},h\}$ as above. Then:
\begin{enumerate}
\item $\mu_{f}(z) = \overline{\mu_{f}(\bar{z})}$.
\item Let $f_{u}$ denote the restriction of $f$ to $\mathbb{H}$. Then $(\pi_{Q})^{-1} \circ f_{u} \circ \pi_{P}$ is the solution to Problem~\ref{polygonmappingproblem} with data $P$ and $Q$.
\end{enumerate}
\end{lemma}

\noindent\textit{Proof:} We first prove that for all $z \in \mathbb{C}$, $f(z) = \overline{f(\bar{z})}$. Define another homeomorphism $g$ as $g(z)=\overline{f(\bar{z})}$. It is straightforward to check that $g$ is a self homeomorphism of $\hat{\mathbb{C}}$ and satisfies $g(z_{i})=w_{i}$.

Now $||\mu_{f}||_{\infty} = ||\mu_{g}||_{\infty}$. By uniqueness of the extremal quasiconformal mapping, $f$ is unique, and so must satisfy $f = g$ everywhere. Thus $f(z) = \overline{f(\bar{z})}$, which implies $\mu_{f}(z) = \overline{\mu_{f}(\bar{z})}$.

To prove the second assertion, let $f^{*}$ denote the solution to Problem~\ref{polygonmappingproblem} with data $P$ and $Q$. Using the above construction, we get a self-homeomorphism $h^{*}$ of $\hat{\mathbb{C}}$ which satisfies the same properties as $h$ and $f$. Uniqueness of $f$ now implies that $f = h^{*}$.
\end{proofof}
}

\paragraph{Ways to represent the Teichm\"{u}ller map}
In theory, a normalized q.c.h. $f$ can be specified by specifying $\mu_{f}$. For computational purposes, unless a closed form expression for $f_*$ or $\mu_*$ is available, the best one can do is to evaluate $f_*$ or $\mu_{*}$ at a dense set of point inside the domain. Teichm\"{u}ller's characterization states that $\mu_*$(the Beltrami coefficient of the solution to either Problem~\ref{polygonmappingproblem} or Problem~\ref{puncturedsphereproblem}), equals $ k | \phi| / \phi$, for some $0 \leq k <1$, and some \textit{integrable holomorphic quadratic differential} $\phi$ (\defref{hqddef}). $\phi$ comes from an $n-3$ dimensional family, and a closed form expression for a basis $\{\phi_{1},\cdots, \phi_{n-3}\}$ \textit{is} available.
Therefore, by representing the coefficients $c_{i}$ in $\phi = \sum_{i=1}^{n-3} c_{i} \phi_{i}$ and $k$, one can represent $\mu_*$. The input and output complexity of both problems would be $O(n)$ in this case. The q.c.h. $f_*$ is the solution of the Beltrami equation for $\mu_{*}$, and can be represented as a series of singular operators applied to $\mu_{*}$ (\cite{Daripa1993355},\cite{doi:10.1137/050640710},\cite{Ahlfors_2006}).

\paragraph{Our representation:} We do not perform a search on the coefficients $c_{i}$ and adopt the first approach instead. If $k_*= ||\mu_{*}||_{\infty}$ is the maximal dilatation of the extremal map $f_*$, then our goal can be stated as follows.

\textbf{Goal}: Given an $\epsilon>0$, compute the values of $f$ on a given set of points inside the base polygon $P$, where the Beltrami differential $\mu_{f}$ of $f$ satisfies $||\mu_{f}||_{\infty} < k_*+\epsilon $.

\paragraph{Complexity:} To the best of the authors' knowledge, even if the polygons $P$ and $Q$ have rational coordinates, there is no known way to represent the extremal map with finite precision (all representations may consist of transcendental numbers). Thus, it is not known whether the problem is in NP or not. We therefore do not address the  actual complexity, and straightaway aim towards an approximation algorithm.
%\mayank{@Vamsi: Do we know how big the size of the mesh should be to get $\epsilon$ accuracy? Even a loose upper bound would do. I would replace dense in the above Goal with that upper bound then.}
%\vamsi{@Mayank: Mate, I am not sure of the size yet. Besides, the size of the input mesh does not really affect the answer. I mean, if you give me a coarse mesh, I will construct a fine one in the algorithm and just return the images of the coarse mesh. Also, we aren't returning the images of these points exactly. Even that is approximate with some error tolerance. (It is difficult to estimate the error caused, but I think (I have to write down to be sure) I can do so in terms of certain constants that still need to be estimated. For instance one of those constants is the constant in the Caldeon-Zygmund inequality.) By the way it is an additive error.}

\subsection{Our results}
\subseclabel{ourresults}

\textbf{Continuous construction:}
Problem~\ref{puncturedsphereproblem} asks for the extremal Beltrami differential on $R$ (the $z_{i}$ punctured sphere) that is isotopic to the starting map $h$. All Beltrami differentials of q.c.h. that are isotopic to $h$ (relative to the punctures) constitute what is called the global equivalence class (\defref{globalequivalence}) of $\mu_{h}$, and our task is to compute the Beltrami differential in this class with the least $L_{\infty}$ norm. Denote the vector space of all Beltrami differentials on the base $z_{i}$-punctured sphere $R$ by $B(R)$, and the unit ball (in  the $L_{\infty}$ norm) of this vector space as $B_{1}(R)$. 

The global class of $\mu_{h}$ cannot be described in a closed form (the only way to know if two differentials are globally equivalent is to solve their Beltrami equation). It lies inside $B_{1}(R)$, and except in trivial cases, is not convex.

Our main result is that we solve the problem by breaking the $L_{\infty}$ minimization over the global class into a sequence of $L_\infty$ minimizations over a convex domain $\mathcal{D}(\mu_{h})$, described explicitly (in terms of $O(n)$ equalities) in terms of $\mu_{h}$. This convex domain will be the class of Beltrami differentials that are \textit{infinitesimally equivalent} (\defref{infinitesimalequivalence}) to $\mu_{h}$. 

Let $\mathcal{D}(\mu)$ denote the infinitesimal equivalence class of $\mu$, and $P(\mu)$ the Beltrami differential $\nu_{0} \in \mathcal{D}(\mu)$ such that $||\nu_{0}||_{\infty} \leq ||\nu||_{\infty}$ for all $\nu \in \mathcal{D}(\mu)$. $\mathcal{P}(\mu)$ is called infinitesimally extremal (\defref{infinitesimallyextremal}).

\begin{theorem}\thmlabel{continuous1}[Limiting procedure for Punctured Sphere Problem]
There exists a sequence of q.c.h. $f_{i}$ s.t.
\begin{CompactEnumerate}
\item $f_{1}=h$, the starting map in Problem~\ref{puncturedsphereproblem}.
\item \textbf{Isotopic:} $f_{i}$ is homotopic to $h$, and$f_{i}(z_{j})=w_{j}$, for all $i$ and  $j$.
\item \textbf{``Explicit" construction:} Denote by $\mu_{i}$ the Beltrami coefficient of $f_{i}$. Then $\mu_{i+1}$ is an ``explicit function" of $\mu_{i}$ and $P(\mu_{i})$ in that it can be obtained by solving two differential equations depending only on $\mu_{i}$ and $P(\mu_{i})$.
\item \textbf{Uniform Convergence:} $f_{i} \rightarrow f_*$ uniformly and $\Vert \mu_i \Vert _{L^{\infty}} \rightarrow \Vert \mu_{*} \Vert _{L^{\infty}}$ as $i \rightarrow \infty$.
\end{CompactEnumerate}

%This is the theorem stating that there is a sequence of q.c.h. converging in the limit to the extremal map. Obviously, we need to state why this sequence is special, i.e., $f_{n+1}$ depends on $f_{n}$ and $\mu_{fn}$ only. \vamsi{Sure, maybe we can say that $f_{n+1} = v_n \circ s_n \circ f_n$ where $s_n$ is the solution of a Beltrami equation whose Beltrami coefficient depends on $\mu_n$ and $v_n$ is a diffeomorphism obtained by finding the trajectories of a vector field dependent on $s_n$. Now $f_{n} \rightarrow f$ uniformly and so does $\Vert \mu_n \Vert _{L^{\infty}} \rightarrow \Vert \mu \Vert _{L^{\infty}}$.}

\end{theorem}

\begin{theorem}\thmlabel{continuous2}[Fast approximation]
There exist constants $C>0$ and $\epsilon_{0}>0$ such that for all $\epsilon < \epsilon_{0}$ and
$$\forall \ \ i \geq \frac{C}{ \epsilon ^{4} (1-||\mu_{1}||_{\infty})^{2}} , $$

$||\mu_{i}||_{\infty} - k_{*} < \epsilon$, where $\mu_{i}$ is the Beltrami differential of $f_{i}$ in \thmref{continuous1} above.
\end{theorem} 
%\mayank{I think $10$ can be replaced by a $2$ here. See Proof of \lemref{one_step_decrease_final}}

\noindent\textbf{Discretization:} 
We represent all Beltrami differentials as piecewise constant differentials\footnote{In fact, the existence of the solution to the Beltrami equation of an arbitrary $\mu \in L^{\infty}$ with $||\mu||_{\infty}<1$ was shown by 1) first showing the existence of the solution to a piecewise constant $\mu^{'}$, 2) sewing the individual piecewise q.c. maps along the boundary, and 3) taking a limit of such piecewise constant differential $\mu_{n}^{'} \rightarrow \mu$ and showing that the maps converge.} on a fine mesh. Every step of the continuous procedure mentioned above is shown to have a discrete analog.

The mesh we will be working on depends on the error tolerance $\delta$ required; near the punctures the mesh is made up of (triangulated) regular polygons, whose number of vertices and radii depend on $\delta$. The mesh is a triangulation with edge lengths bounded above by $\epsilon$ (which is a function of $\delta$). We call this triangulation a canonical triangulation of size $\epsilon$ (see \defref{canonicaltriangulation}) and denote it by $\Delta_{\epsilon}$.

% The following algorithms run on an analogue computer that can compute sums, products, quotients, and subtractions of arbitrary real numbers at unit cost.
The first theorem states that our discretization for the operator $\mathcal{P}$ that returns the infinitesimally extremal Beltrami coefficient is in fact an approximation.

\begin{theorem} \thmlabel{discrete1}[Discrete infinitesimally extremal]
Given an error tolerance $0<\delta<1$, a collection of $n$ punctures $z_1, z_2, \ldots z_n$, a triangulation $\Delta_{\epsilon}$ and a piecewise constant Beltrami coefficient $\mu$ (where $\Vert \mu \Vert <1$), there exists an algorithm \INEXT that computes a piecewise constant Beltrami coefficient $\hat{\nu}$ such that $\Vert \hat{\nu} \Vert - \Vert \nu \Vert <  \Vert \nu \Vert \delta$, where $\nu = \mathcal{P}(\mu)$.
\end{theorem}

\noindent\textbf{Discrete algorithm:}
Having discretized the main component of our procedure, all the other steps in our procedure can be easily implemented in practice. Computational quasiconformal theory is a field still in its infancy, and very few error estimates on these widely-used discretizations are known. We define two subroutines next that concern the discretization of compositions and inverses of quasiconformal maps. 

\begin{definition}\deflabel{PIECEWISECOMP}[Subroutine: \PIECEWISECOMP]

\noindent \textbf{Input:}  A triangulation $\Delta_{\epsilon}$, two piece-wise constant Beltrami coefficients $\mu_1$ and $\mu_2$ (corresponding to q.c.h $f_1$ and $f_2$ respectively), and error tolerances $\delta_1$ and $\delta_2$.

\noindent \textbf{Output:} A triangulation $\Delta_{\epsilon^{'}}$ that is a refinement of $\Delta_{\epsilon}$, a piecewise constant Beltrami coefficient $\mu_{\text{comp}}$ that approximates the Beltrami coefficient of the composition $f_3 = f_1 \circ f_2$ within error $\delta_1$ in the $L^{\infty}$ topology, and the images $f_3 (v_a)$ of the vertices $v_a$ of $\Delta_{\epsilon ^{'}}$ up to an error of $\delta_2$.
\end{definition}

\begin{definition}\deflabel{PIECEWISEINV}[Subroutine: \PIECEWISEINV]

\noindent \textbf{Input:}  A triangulation $\Delta_{\epsilon}$, a piecewise constant Beltrami coefficient $\mu$ (corresponding to q.c.h $f$), and error tolerances $\delta_1$ and $\delta_2$.

\noindent \textbf{Output:} A triangulation $\Delta_{\epsilon^{'}}$ that is a refinement of $\Delta_{\epsilon}$, a piecewise constant Beltrami coefficient $\mu_{\text{inv}}$ that approximates the Beltrami coefficient of $f^{-1}$ within error $\delta_1$ in the $L^{\infty}$ topology, and the images $f^{-1}(v_a)$ of the vertices of $\Delta_{\epsilon^{'}}$ up to an error of $\delta_2$.
\end{definition}

Assuming the existence of the subroutines \PIECEWISECOMP and \PIECEWISEINV, we construct an approximation algorithm for the Teichm\"{u}ller map.

\begin{theorem} \thmlabel{discrete2}[Teichm\"{u}ller Map Algorithm]
Given 1) a triangulation $T_0$ that includes $n$ punctures $z_1,\ldots z_n$, 2) a mesh of sample points $S$, 3) an error tolerance $\delta$, and 4) a piece-wise constant Beltrami coefficient $\mu_0$, whose corresponding q.c.h. $f_0$ satisfies $f_{0}(z_{j}) = w_{j}$,
there exists an algorithm \EXTREMAL that computes $\Delta_{\epsilon}$, and the images of $S$ up to an error of  $\delta$ under a q.c.h. $F$ having a piece-wise constant (in the computed triangulation) Beltrami coefficient $\mu_F$ such that 
\begin{CompactEnumerate}
\item $\Vert \mu _F \Vert _{L^{\infty}} - \Vert \mu_{*} \Vert _{L^{\infty}} < \delta$ where $\mu _ *$ is the Beltrami coefficient of the extremal quasiconformal map on the punctured sphere in the homotopy class of $f_0$. 
\item $|F(z_{i})-w_{i}| = O(\delta)$.
\end{CompactEnumerate}
\end{theorem}

An implementation of our algorithm will be presented in a forthcoming paper.

\textbf{Structure of the paper:} In \secref{prelim} we define some terms that we use in our construction. In \secref{continuous} we dwelve into the proofs of \thmreftwo{continuous1}{continuous2}. \secref{discrete} describes our discretized procedure and proves \thmreftwo{discrete1}{discrete2}. We conclude in \secref{discussion} with discussions on complexity and generalizations to arbitrary Riemann surfaces.

\section{Preliminaries}
\seclabel{prelim}

\subsecref{rmtrs} and \subsecref{qcmapsapp} (in Appendix \secref{app1}) provide the basic definitions of Riemann surfaces and quasiconformal maps, respectively. For the sake of completeness of the main body, in this section we define some of the concepts we will require for our continuous construction.

\begin{definition}\deflabel{hqddef}[Holomorphic quadratic differential]
A holomorphic quadratic differential on a Riemann surface $R$ is an assignment of a function $\phi_{i}(z_{i})$ on each chart $z_{i}$ such that if $z_j$ is another local coordinate, then $\phi_{i}(z_{i})=\phi_{j}(z_{j}){(\frac{dz_{j}}{dz_{i}})}^{2}$.
\end{definition}

We will denote the space of such differentials on $R$ as $A(R)$. By the Riemann-Roch theorem, the complex dimension of this  vector space for a genus $g$ closed compact surface with $n$ punctures is $3g-3+n$.

\textbf{Fact:} For $R = \hat{\mathbb{C}} \setminus \{0,1,\infty,z_{1},...z_{n-3}\}$ (the Riemann sphere with $n$ punctures),
\begin{equation} \eqlabel{puncturediff}
\phi_{k}(z) = \frac{\eta_{i}}{z(z-1)(z-z_{k})}, \ \ \ 1\leq k \leq n-3,
\end{equation}
form a basis of $(n-3)$ dimensional complex vector space $A(R)$. Here $\eta_{i}$ is a constant, chosen such that the norm of $\phi$ is $1$, i.e., $||\phi||= \int_{R} |\phi| = 1$. 

\paragraph{Equivalence relations on Beltrami coefficients}

\noindent\textbf{Global equivalence:} This relation is defined only on Beltrami differentials of norm less than $1$, i.e. those that belong to the unit ball $B_{1}(R) = \{ \mu \in B(R) : ||\mu||_{\infty} <1 \}.$ Given two such differentials $\mu$ and $\upsilon$, denote the solution to their respective normalized\footnote{Fixing the points $0$,$1$ and $\infty$. Hence the freedom of M\"{o}bius tranformation is accounted for.} Beltrami equations as $f^{\mu}: R \rightarrow R_{0}$ and $f^{\upsilon}: R \rightarrow R_{1}$. Both $R_{0}$ and $R_{1}$ are punctured spheres.

\begin{definition}\deflabel{globalequivalence}[Global equivalence]
$\mu$ and $\upsilon$ are called globally equivalent (written $\mu \sim_{g} \upsilon$) if:
\begin{CompactEnumerate}
\item $f^{\mu}(z_{i}) = f^{\upsilon}(z_{i})$  $\forall i$.
\item The identity map from $R_{0}$ to $R_{1}$ is homotopic to $f^{\upsilon} \circ (f^{\mu})^{-1}$ via a homotopy consisting of quasiconformal homeomorphisms.
\end{CompactEnumerate}
\end{definition}

\begin{definition}\deflabel{trivial}[Trivial Beltrami differential]
A Beltrami differential $\upsilon$ is called trivial if it is globally equivalent to $0$.
\end{definition}

\noindent\textbf{Infinitesimal Equivalence:} This relation is defined on all of $B(R)$.
\begin{definition}\deflabel{infinitesimalequivalence}[Infinitesimal equivalence]
$\mu$ and $\upsilon$ are infinitesimally equivalent (written $\mu \sim_{i} \upsilon$) if $\int_{R} \mu \phi = \int_{R} \upsilon \phi$ for all $\phi \in A(R)$, with $||\phi||= 1$.
\end{definition}

\begin{definition}\deflabel{infinitesimallyextremal}[Infinitesimally extremal]
A Beltrami differential $\upsilon$ is called infinitesimally extremal if $||\upsilon||_{\infty} \leq ||\mu||_{\infty}$ for all $\mu \sim_{i} \upsilon$.
\end{definition}

\begin{definition}\deflabel{infinitesimallytrivial}[Infinitesimally trivial Beltrami differential]
A Beltrami differential $\upsilon$ is called infinitesimally trivial if it is infinitesimally equivalent to $0$.
\end{definition}

\subsecref{ingredients} lists all the theorems from Teichm\"{u}ller theory that we will require, namely the Mapping theorem (\thmref{mapping}), the composition formula (Equation~\eqreff{formula}), the variational lemma (\lemref{variationallemma}), the principle of Teichm\"{u}ller contraction (Equation~\eqreff{tcontract}) and most importantly, the Hamilton-Krushkal, Reich-Strebel, necessary and sufficient condition for optimality (\thmref{hkrs}).

\section{Proofs of \thmreftwo{continuous1}{continuous2}}
\seclabel{continuous}

At the heart of our construction of the sequence $f_{i}$ in \thmreftwo{continuous1}{continuous2} lies the following lemma. Let $h$ be any quasiconformal homeomorphism between $R$ (the $z_{i}$ punctured sphere) and $S$ (the $w_{i}$ punctured sphere) which is a valid input to Problem~\ref{puncturedsphereproblem}, and $\mu_{h}$ denote its Beltrami differential.

\begin{lemma}\lemlabel{reductionlemma}
Let $\upsilon_{h}$ be the infinitesimally extremal Beltrami differential in the infinitesimal class of $\mu_{h}$. Let $\mu_{g}(t)$ be a curve of Beltrami differentials with the following properties:

\begin{CompactEnumerate}
\item $\mu_{g}(t)$ is globally trivial.
\item $\mu_{g}(t) = t(\mu_{h} - \upsilon_{h}) + O(t^{2})$.
\end{CompactEnumerate}

Denote the solution to the Beltrami equation of $\mu_{g}(t)$ by $g_{t}$. Then $\exists \ \delta>0$ such that $\forall t <\delta$, the map $h_{t} = h \circ (g_{t})^{-1}$ has smaller dilatation than $h$.
\end{lemma}

\textlater{\begin{proofof}{\lemref{reductionlemma}}
By the formula for composition of quasiconformal maps (\eqreff{formula} in \subsecref{ingredients}),

\begin{equation}\eqlabel{equ1}
 \mu_{h_{t}}(g_{t}(z)) = \frac{\mu_{h}-\mu_{g}(t)}{1-\overline{\mu_{g}(t)}\mu_{h}}\frac{1}{\theta_{t}},
\end{equation}

where $\theta_{t} = \overline{p_{t}}/p_{t}$ and $p_{t} = \frac{\partial g_{t}}{\partial z}$. \eqreff{equ1} implies

\begin{equation} \eqlabel{equ0}
| \mu_{h_{t}} \circ g_{t} |^{2} = \frac{|\mu_{h}|^{2} - 2 \text{Re} (\mu_{h} \overline{\mu_{g}(t)}) +|\mu_{g}(t)|^{2} }  { 1 - 2 \text{Re}( \mu_{h} \overline{\mu_{g}(t)}) +|\mu_{g}(t) \mu_{h}|^{2}                }
\end{equation}

Using the fact that
\[|| \mu_{g}(t) -t (\mu_{h} - \upsilon_{h})||_{\infty} = O(t^{2}) \]
and differentiating \eqreff{equ0} with respect to $t$ once and putting $t=0$, we get that

\begin{equation}\eqlabel{equ2}
|\mu_{h_{t}} \circ g_{t}| = |\mu_{h}| -t  \frac{1-|\mu_{h}|^{2}}{|\mu_{h}|} \text{Re}(|\mu_{h}|^{2} - \mu_{h} \overline{\upsilon_{h}}) + O(t^{2})
\end{equation}

Let $k_{0} = ||\upsilon_{h}||_{\infty} < k = ||\mu_{h}||_{\infty}$. Define

\[ S_{1} = \{ z \in R : |\mu_{h}(z)| \leq (k+k_{0})/2 \} \] and
\[ S_{2} = \{ z \in R : (k+k_{0})/2 < |\mu_{h}(z)| \leq k \} \]

Clearly, $S_{1} \cup S_{2} = R$. Since in $S_{1}$ the starting value of this curve at $t=0$ is $|\mu|$, which is certainly less than $k$, \eqreff{equ1} implies there exists $\delta_{1} >0$ and $c_{1}>0$ such that for $0<t <\delta_{1}$,

\begin{equation}\eqlabel{equ3}
 |\mu_{h_{t}} \circ g_{t} (z)| \leq k - c_{1}t \ \ \ \text{for} \ z \in S_{1}
\end{equation}

For $z$ in $S_{2}$ the coefficient of $t$ in \eqreff{equ2} is bounded below by

\[ \frac{1-k^{2}}{k} \left[ {\left( \frac{k+k_{0}}{2} \right)}^{2} -k_{0}k \right] = \frac{1-k^{2}}{k} {\left( \frac{k-k_{0}}{2} \right)}^{2} >0 \]

Therefore, \eqreff{equ2} implies there exists $\delta_{2} >0$ and $c_{2} >0$ such that for $0 <t < \delta_{2}$,

\begin{equation}\eqlabel{equ4}
 |\mu_{h_{t}} \circ g_{t} (z)| \leq k - c_{2}t \ \ \ \text{for} \ z \in S_{2}
\end{equation}

Putting together \eqreff{equ3} and \eqreff{equ4}, we find that $||\mu_{h_{t}}||_{\infty} < k$ for sufficiently small $t>0$, proving the lemma.

\end{proofof}
}

The proof is similar to that of the Hamilton-Krushkal, Reich-Strebel necessary-and-sufficient condition for extremality (see \thmref{hkrs}), published in a sequence of papers. We refer the reader to \cite{Gardiner} for a combined proof of this celebrated result, which is the one we adapt. To the best of the authors' knowledge, the above lemma is the first result that describes, given a starting map, how to get a map with a smaller dilatation.

The proof of \thmref{continuous1} is constructive. We summarize the construction first:

\subsection{Summary of the construction}\subseclabel{summarycontinuous}
At step $i$, Given a starting map $f_{i}:R \rightarrow S$ with Beltrami coefficient $\mu_{i}$, let $\upsilon_{i}$ denote the infinitesimally extremal Beltrami coefficient in the infinitesimal class of $\mu_{i}$. Let $k_{i} = ||\mu_{i}||_{\infty}$ and $k_{i}^{0}=||\upsilon_{i}||_{\infty}$. Observe that $\mu_{i} - \upsilon_{i}$ is infinitesimally trivial (\defref{infinitesimallytrivial}).

\begin{CompactEnumerate}
\item Choose $t$ such that
\begin{equation}\eqlabel{choosing_t}
t =\min \left( \frac{3}{4},C_{1},\frac{\epsilon}{4},\sqrt{\frac{\epsilon}{2C_{2}}},\frac{(k_{i}-k_{i}^{0})^2(1-k_{i}^2)}{1-k_{i}^2+C_{2}} \right),
\end{equation}
where $\epsilon \leq \min (1/2, (k_{i}-k_{i}^{0})/8)$, and $C_{1}$ and $C_{2}$ are two explicit constants to be derived later.

\item Use \subsecref{vectorfieldcontinuous} to construct a quasiconformal self-homeomorphism $g_{i}$ of $R$ such that 
\begin{CompactItemize}
\item $\mu_{g}$ is globally trivial.
\item $||\mu_{g}-t(\mu_{i}-\upsilon_{i})||_{\infty}< C_{2} t^{2}$, where $C_{2}$ is the same constant as in \eqreff{choosing_t}.
\end{CompactItemize}
\item Form $f_{i+1} = f_{i} \circ (g_{i})^{-1}$ such that $f_{i+1}$ has smaller dilatation than $f_{i}$ (by \lemref{reductionlemma}).
\item Reiterate with $f_{i+1}$ as the starting map.
\end{CompactEnumerate}

\subsection{How the construction implies \thmreftwo{continuous1}{continuous2}}
\subseclabel{proofs_of_theorems}
The $\mu_{g}$ in step (2) above can be constructed by solving two differential equations involving $t$, $\mu_{i}$ and $\upsilon_{i}$ (\lemref{mugtheorem} in \subsecref{vectorfieldcontinuous}). Assuming that, we have the following lemma that quantifies the progress made in step $i$. Recall that $k_{i} = ||\mu_{i}||_{\infty}$ and $k_{i}^{0} = ||\upsilon_{i}||_{\infty}$.

\begin{lemma}\lemlabel{one_step_decrease}[Decrease in one step]
If $t$ is chosen as in \eqreff{choosing_t}, then $k_{i}-k_{i+1} > d,$ where

\begin{equation}\eqlabel{one_step_decrease}
d = \min \left(\frac{k_{i}-k_{i}^{0}}{4}, \frac{(k_{i}-k_{i}^{0})^2 t (1-k_{i}^2)}{8} \right). \notag
\end{equation}

\end{lemma}

\textlater{\begin{proofof}{\lemref{one_step_decrease}}
To simplify notation, we let $k = k_{i}$ and $k_{0}= k_{i}^{0}$, since we will be assuming that we are in step $i$ of the iteration. As in the proof of \lemref{reductionlemma}, let $S_1$ be the region where $\vert \mu_{i} \vert \leq \frac{k+k_0}{2}$ and $S_2$ be such that $\frac{k+k_0}{2} \leq \vert \mu_i \vert < k$. Assuming that $t < \min (3/4,C_{1})$ implies \lemref{mugtheorem}, so we assume this condition on $t$.

Furthermore, on $S_1$, if $t < \min \left(3/4,C_{1},\frac{\epsilon}{4}, \sqrt{\frac{{\epsilon}}{{2C_{2}}}} \right)$, by the composition formula for q.c. maps we get,

\begin{eqnarray}
\vert \mu_{i+1} \circ g_t (z) \vert &=& \frac{\vert \mu_{i} - \mu_{g} \vert}{\vert 1 -\mu_i \bar{\mu _g} \vert} \nonumber \\
										&\leq& \frac{\vert \mu _i - t(\mu_i-\nu_i)\vert}{\vert 1 -\mu_i \bar{\mu _g} \vert}+\frac{C_{2}t^2}{\vert 1 -\mu_i \bar{\mu _g} \vert} \nonumber \\
&\leq& \frac{1}{1-\epsilon}\left ( \frac{k+k_0}{2} + 2t + C_{2}t^2 \right ) \nonumber \\
&\leq& \frac{1}{1-\epsilon}\left ( \frac{k+k_0}{2} + \epsilon \right ) \nonumber
\end{eqnarray}
where the last inequality follows by requiring $\vert \mu_i \bar{\mu_{g}} \vert < 2t + C_{2}t^2 < \epsilon$, which is true for the assumed value of $t$. Notice that $\vert \mu_{i} \vert, \vert \nu_{i} \vert $ are less than $1$.

Therefore, on $S_1$,
\begin{eqnarray}
k-\vert \mu _{i+1} \vert &\geq& \frac{k-k_0-2\epsilon (1+k)}{2(1-\epsilon)} \nonumber \\
&>& \frac{k-k_0}{4}\nonumber
\end{eqnarray}
if $\epsilon \leq  \frac{k-k_0}{8}$. On $S_2$,
\begin{eqnarray}
\frac{\vert \mu_{i} - \mu_{g} \vert}{\vert 1 -\mu_i \bar{\mu_{g}} \vert} &\leq& \frac{\vert \mu_{i} - t(\mu_i-\nu_i)\vert}{\vert 1 -\mu_i \bar{\mu_g} \vert}+\frac{C_{2}t^2}{\vert 1 -\mu_i \bar{\mu_g} \vert} \nonumber \\
&\leq& \frac{\vert \mu_{i} - t(\mu_i-\nu_i)\vert}{\vert 1 -\mu_i t (\bar{\mu_{i}}-\bar{\nu_{i}} ) \vert - \vert \mu_{i} \vert C_{2}t^2}+\frac{C_{2}t^2}{\vert 1 -\mu_{i} \bar{\mu_g} \vert}\nonumber
\end{eqnarray}
Now,
\begin{eqnarray}
k-\vert \mu_{i+1} \vert &\geq& \vert \mu_{i} \vert - \vert \mu_{i+1}\vert \nonumber \\
&\geq& \frac{\vert \mu_i \vert (\vert 1 -\mu_i t (\bar{\mu_i}-\bar{\nu_i} ) \vert - \vert \mu _i \vert C_{2} t^2) - \vert \mu_i - t(\mu_i - \nu_i) \vert  }{\vert 1 -\mu_i t (\bar{\mu _i}-\bar{\nu_i} ) \vert - \vert \mu_i \vert C_{2} t^2} - \frac{C_{2}t^2}{1-\epsilon} \nonumber \\
&\geq& \frac{A-B}{\vert 1 -\mu_i t (\bar{\mu_i}-\bar{\nu_i} ) \vert - \vert \mu_i \vert C_{2} t^2} - \frac{C_{2}t^2}{1-\epsilon}
\label{inter}
\end{eqnarray}

where $A = \vert \mu_{i} \vert (\vert 1 -\mu_{i} t (\bar{\mu_{i}}-\bar{\nu_i} ) \vert) $ and $B= \vert \mu_{i} \vert C_{2} t^2 - \vert \mu_i - t(\mu_i - \nu_iu) \vert $.
Using
\begin{gather}
A-B = \frac{A^2-B^2}{A+B} \nonumber \\
A+B < 4 \nonumber \\
A^2-B^2 = (1-\vert \mu_i \vert ^2)(2\Re (t \mu_{i} (\bar{\mu_i} - \bar{\nu_i})) - t^2 \vert \mu_i \vert^2 \vert \mu_i - \nu _i \vert^2(1+\vert \mu_i \vert^2)) \nonumber  \\
\geq  (1-\vert \mu_i \vert ^2)(t(k-k_0)^2 - t^2 \vert \mu_i \vert^2 \vert \mu_i - \nu_i \vert^2(1+\vert \mu_i \vert^2))
\label{next}
\end{gather}
Using Equation~\ref{next} in Equation~\ref{inter} and the fact that $\epsilon < 1/2$ we see that
\begin{gather}
k-\vert \mu_{i+1} \vert \geq \frac{1-k^2}{4}  ( t(k-k_0)^2 -8t^2 ) - 2C_{2}t^2 \nonumber \\
\geq \frac{(k-k_0)^2 t (1-k^2)}{8}
\label{fin}
\end{gather}
with the last equation holding if $t < \frac{(k-k_0)^2(1-k^2)}{1-k^2+C_{2}} $, concluding the proof.

\end{proofof}
}

\noindent Now we apply the principle of Teichm\"{u}ller contraction, which essentially bounds $k_{i}-k_{i}^{0}$ from below by a function of $k_{i}-k_{*}$ (how far we are from the infinitesimally extremal coefficient tells us how far we are from the extremal). Using \lemref{one_step_decrease}, Equation~\eqreff{choosing_t} and the principle of Teichm\"{u}ller contraction (Equation~\eqreff{tcontract}), we get

\begin{lemma}\lemlabel{one_step_decrease_final}
There exists a constant $C_{3}>0$, such that if $t$ is chosen as in \eqreff{choosing_t}, then $$k_{i}-k_{i+1} > C_{3}(1-||\mu_{1}||_{\infty})^{2}(k_{i}-k_{*})^{4}$$
\end{lemma}

\textlater{\begin{proofof}{\lemref{one_step_decrease_final}}
The first three terms in the $\min$ expression in Equation~\eqreff{choosing_t} are independent of the iteration step. If the value of $t$ is one of these, then the lemma is evident. Similarly, if the value of $d$ is the first of the two terms in the $\min$ in \lemref{one_step_decrease}, then the lemma is clear too.

Assume now that $t = \frac{(k_{i}-k_{i}^{0})^2(1-k_{i}^2)}{1-k_{i}^2+C_{2}}$. Noting that $1-k_{i}+C_{2}< C_{2}+1$, and that $1-k_{i}^{2}>1-k_{i}$, we get $t > \frac{(k_{i}-k_{i}^{0})^2(1-k_{i})}{1+C_{2}}$. Using this value in \lemref{one_step_decrease} gives

\begin{equation}
d  > \frac{(k_{i}-k_{i}^{0})^{4} (1-k_{i})^{2}}{8(1+C_{2})}, \notag
\end{equation}

and by Teichm\"{u}ller contraction, $k_{i}-k_{i}^{0} \geq (k_{i} - k_{*})/10$, implying
\begin{equation}
d  > \frac{(k_{i}-k_{*})^{4} (1-k_{i})^{2}}{80(1+C_{2})}, \notag
\end{equation}

Putting $C_{3} = \frac{1}{80(1+C_{2})}$, and noting that $1-k_{i} > 1- k_{1}$ completes the proof of the lemma.

\end{proofof}
}

Using \lemref{one_step_decrease_final}, the proofs of \thmreftwo{continuous1}{continuous2} can now be completed.

\noindent\textbf{Proof of \thmref{continuous1}:} The first three assertions follow easily from our construction (notice that $\mu_{g}$ being trivial implies that we stay within the initial homotopy class). We now prove assertion (4) in the statement of the theorem.

The $k_{i}$ form a decreasing sequence and are bounded from above by $k_{1}$ and below by $k_{*}$. Hence the $k_{i}$ converge to some $k\geq0$. But then $k_{i}-k_{i+1} \rightarrow 0$ as $i \rightarrow \infty$, and \lemref{one_step_decrease_final} now implies that $k_{i}-k_{*} \rightarrow 0$, and so $k=k_{*}$. Thus $k_{i} \rightarrow  k_{*}$.

Now the fact that $f_{i}$ converge uniformly to $f_{*}$ follows because of the convergence property of q.c.h. The space of all q.c.h. $f$ with $||\mu_{f}||_{\infty} < k_{1}$ forms a compact space, and so there exists a subsequence that converges. By the arguments above and by uniqueness of $f_{*}$, this limit must be $f_{*}$. Furthermore, this must be true of any convergent subsequence of the $f_{i}$. Thus we get that the entire sequence $f_{i}$ converges uniformly to $f_{*}$.

\noindent\textbf{Proof of \thmref{continuous2}:} 
Let $A= C_{3}(1-\Vert \mu_1 \Vert _{\infty}) ^{2}$ and define $y_n = \Vert \mu_n \Vert _{\infty} - k_{*}$. By \lemref{one_step_decrease_final} and Teichm\"{u}ller contraction (Equation~\eqreff{tcontract}), $y_{n+1} \leq y_n - A y_n ^4$ . If $y_1 <\epsilon$, we are done. If not, then $y_1 - y_2 \geq A y_1 ^4 \geq A \epsilon ^4$ and thus $y_2 \leq y_1 - A \epsilon ^4$. If $y_2 < \epsilon$ we are done. If not, $y_3 \leq y_2 - A \epsilon ^4 \leq y_1 - 2A \epsilon ^4$. Continuing inductively we see that $y_n \leq y_1 - (n-1) A \epsilon ^4$ if $y_{n-1} > \epsilon$. The right hand side is less than $\epsilon$ if $n>\frac{1}{A \epsilon ^4} (\Vert \mu _ 1 \Vert _{L^{\infty}} - k_{*} - \epsilon)$. Since $(\Vert \mu _ 1 \Vert _{L^{\infty}} - k_{*} - \epsilon) <2$, putting $C = 2/C_{3}$ proves the theorem.

\subsection{Constructing self homeomorphisms $g_{i}$}
\subseclabel{vectorfieldcontinuous}
Given a starting map $f_{i}$, we show how to construct the self homeomorphism $g_{i}$ of $R$ used in our construction. We simplify notation by suppressing the $i$, keeping in mind that this is the $i$th step of the procedure. Thus $\mu$ and $\mu_{g}$ will denote the Beltrami differentials of $f_{i}$ and $g_{i}$, respectively. Also, $\upsilon$ is the infinitesimally extremal Beltrami differential in the infinitesimal class of $\mu$.

Let $\alpha= \mu-\upsilon$, $t$ be as in Equation~\eqreff{choosing_t}, and let $f^{t\alpha}$ be the normalized solution to the Beltrami equation for $t \alpha$. We observe next that $f^{t\alpha}$ moves the points $z_{i}$ only by a distance $O(t^{2})$.

\begin{lemma}\lemlabel{shift}
Let $r=\displaystyle \max_{1 \leq i \leq n-3} |z_{i}|$, and let $f^{t \alpha}$ be as above. Then there exists a constant $C_{r}$ depending only on $r$, and a constant $\delta>0$ such that for all $i$, $|f^{t \alpha}(z_{i}) - z_{i}| \leq C_{r} t^{2}, \ \ \forall t<\delta$.
\end{lemma}

\textlater{\begin{proofof}{\lemref{shift}}
By \eqreff{puncturediff}, \[\phi_{i}(\zeta) = \frac{z_{k}(z_{k}-1)}{\zeta(\zeta-1)(\zeta-z_{k})}\] for $1 \leq i \leq n-3$ is a basis for the space of quadratic differentials on $R$. Let $\zeta=\xi + i \eta$. Infinitesimal equivalence of $t\alpha$ now implies that
\begin{equation}\eqlabel{pairingiszero}
\int \int_{\mathbb{C}} \frac{t \alpha d \xi d \eta}{\zeta(\zeta-1)(\zeta-z_{k})} = 0
\end{equation}
Now we use the mapping theorem (\thmref{mapping} in \subsecref{ingredients}). In the notation of the theorem, $V(z_{i})=0$ by \eqreff{pairingiszero}. Existence of $\delta$, $C_{r}$ and the statement now follows from the statement of the mapping theorem.
\end{proofof}
}
Denote $f^{t \alpha}(z_{i})$ by $z_{i}^{'}$. We will first construct another homeomorphism $f_{\textbf{v}}$ from $\hat{\mathbb{C}}$ to itself which satisfies $f_{\textbf{v}}(z_{i}^{'}) = z_{i}$. We then define the required self homeomorphism $g_{i}=f_{\textbf{v}} \circ f^{t \alpha}$. The construction of $f_{\textbf{v}}$ will be via a vector field method.

\paragraph{Construction of $f_{\textbf{v}}$ by a vector field method:}
Let $\{D_{1},\cdots,D_{n-3}\}$ denote disjoint open disks centered at $z_{i}$. Choosing the radius of each disk to be $r=d/4$, where $d= \displaystyle \max_{1 \leq i,j \leq n-3} |z_{i}-z_{j}|$ ensures disjointness. We will fix these disks once and for all.

\paragraph{A single disk:}We first construct a self homeomorphism $f_{\textbf{v}}^{i}$ of $\hat{\mathbb{C}}$ which is the identity outside $D_{i}$, and maps $z_{i}^{'}$ to $z_{i}$. Now $z_{i} \in \mathbb{R}$, and by a rotation we can assume that $z_{i}^{'}$ is real and greater than $z_{i}$. Consider the vector field \[ X(z) = p(z) (z_{i}^{'}-z_{i})\frac{\partial }{\partial x}, \]
where $p(z)$ is a $C^{\infty}$ function identically zero outside $D_{i}$, and identically $1$ inside the disk of radius $r/2$ around $z_{i}$, denoted as $D_{i}^{'}$.
Let $\gamma$ be the one parameter family of diffeomorphisms associated with this vector field. We denote the time parameter by $s$ and note that the diffeomorphism $\gamma_{1}$ sends $z_{i}^{'}$ to $z_{i}$. We denote this diffeomorphism $\gamma$ at $s=1$ by $f_{\textbf{v}}^{i}$. Now define $f_{\textbf{v}} = f_{\textbf{v}}^{n-3} \circ f_{\textbf{v}}^{n-2} \cdots \circ f_{\textbf{v}}^{1}$, and $g_{i}=f_{\textbf{v}} \circ f^{t \alpha}$. We have

\begin{lemma}\lemlabel{mugtheorem}
There exist constants $C_{1}>0$ and $C_{2}>0$ such that $\mu_{g}$ above is globally trivial, and  for all $t < \min(3/4, C_{1})$, $||\mu_{g}-t(\mu-\upsilon)||_{\infty}< C_{2} t^{2}.$
\end{lemma}
The exact values of $C_{1}$ and $C_{2}$ can be inferred from the proof. They equal the values of $\delta$ and $C_r$, respectively, in the mapping theorem, when $r= 1$.

\textlater{\begin{proofof}{\lemref{mugtheorem}}
%\documentclass[a4paper]{article}
%\usepackage{amssymb}    %use with LaTeX2e
%\usepackage{amsthm}
%\usepackage{fullpage}
%\newtheorem{thm}{Theorem}
%\newtheorem{proposition}{Proposition}[section]
%\newtheorem{lemma}{Lemma}[section]
%\begin{document}

%\label{Proving the existence of C and finding it}

\indent In what follows, we denote $f_i$ (the quasiconformal map after $i$ iterations of the algorithm) by $f$ for convenience. In fact, we drop the subscript $i$ altogether. Recall that from $\mu_i$, we construct a new quadratic differential $\nu_i$ that is infinitesimally extremal. Let $\alpha _i = \mu_i - \nu_i$.\\
\indent  We then construct $g_i$ that is a self homeomorphism of the base punctured sphere $R$ using the vector field method. If we manage to prove that the dilatation of the map obtained by the vector field method $\Vert v \Vert _{L^{\infty}} \leq G t^2$, then for the composition (here $\vert A \vert =1$) $g_i = f_{\mathbf{v}} \circ f^{t\alpha _i}$, if $i$ is sufficiently large then $\vert 1+t\alpha v A \vert >\frac{1}{2}$  
\begin{eqnarray}
\vert \frac{t\alpha + A v \circ f^{t\alpha}}{1+t \alpha A v} - t\alpha \vert &=& \vert \frac{Av \circ f^{t\alpha} + t^2 \alpha ^2 A v \circ f^{t\alpha}}{1+t\alpha v A}\vert  \nonumber \\
&\leq& (2G+1)t^2 
\label{stringof}
\end{eqnarray}
\indent Therefore the Beltrami of the composition is $t\alpha + O(t^2)$ where the  $O(t^2)$ term is bounded above by $(2G+1)t^2=Ct^2$. Our aim is to prove that $G$ exists and is bounded independent of $i$. In whatever follows we denote $f_{\mathbf{v}}$ simply by $f$. \\
\indent Recall that the vector field is $X = \sum_j \rho_j (z_j (t) - z_j)$. Let $\gamma (t,s,y,\bar{y})$ be the flow where $s$ is the ``time parameter" for the flow and $y$ is the initial position. Notice that the vector field diffeo is $\gamma(t,1,y,\bar{y})=f$ and the Beltrami $v$. Notice that 
\begin{eqnarray}
v &=& \frac{f_{\bar{y}}}{f_y} \nonumber \\
v_t &=& \frac{f_{\bar{y}t}f_{y} - f_{\bar{y}}f_{ty}}{f_y ^2} \nonumber \\
v_{tt} &=& \frac{(f_{\bar{y}tt}f_{y} - f_{\bar{y}}f_{tty})f_y ^2 - 2f_y f_{yt}(f_{\bar{y}t}f_{y} - f_{\bar{y}}f_{ty})}{f_y ^4} \nonumber
\end{eqnarray}
\indent If $\vert f_y \vert > m$ and $f_t$, $\vert f_{y} \vert$, $\vert f_{\bar{y}} \vert$, $\vert f_{\bar{y}tt} \vert$, $\vert f_{ytt} \vert $, $\vert f_{\bar{y}t}\vert$ and $\vert f_{yt} \vert$ are bounded above by $M$, then $\vert v_{tt} \vert \leq \frac{6M^4}{m^4}=G$. First we have to prove that indeed $\vert v \vert < 1$ so that $f$ is q.c. For all of these things we consider
\begin{eqnarray}
\frac{d\gamma}{ds} &=& \sum \rho  (\gamma) (z_i (t)-z_i) \nonumber \\
\frac{d\gamma _t}{ds} &=& \sum (\rho _z  \gamma _t  + \rho _{\bar{z}} \bar{\gamma}_t) (z_i (t)-z_i) + \sum \rho z_i ^{'} \nonumber \\
\frac{d\gamma _{tt}}{ds} &=& \sum (\rho _{zz}  (\gamma _t)^2+2\rho_{\bar{z}z} \vert \gamma _t \vert^2 + \rho_z \gamma _{tt}  + \rho _{\bar{zz}} (\bar{\gamma}_t)^2 + \rho _{\bar{z}} \bar{\gamma} _{tt})  (z_i (t)-z_i) + \sum \rho z_i ^{''} + \sum (\rho _z \gamma _t + \rho _{\bar{z}} \bar{\gamma} _t)z_i ^{'} \nonumber \\
\frac{d\gamma_y}{ds} &=& \sum (\rho _z  \gamma _y  + \rho _{\bar{z}} \bar{\gamma}_y) (z_i (t)-z_i) \nonumber \\
\frac{d\gamma_{\bar{y}}}{ds} &=& \sum (\rho _z  \gamma _{\bar{y}}  + \rho _{\bar{z}} \bar{\gamma}_{\bar{y}}) (z_i (t)-z_i) \nonumber \\
\frac{d\gamma_{yt}}{ds} &=& \sum (\rho _{zz}  \gamma _t \gamma _y + \rho _{\bar{z}z}  \bar{\gamma} _t \gamma _y  + \rho _z  \gamma _{yt}+ \rho _{z\bar{z}}\gamma _t \bar{\gamma}_y + \rho _{\bar{zz}}\bar{\gamma} _t \bar{\gamma}_y + \rho _{\bar{z}} \bar{\gamma}_{yt}) (z_i (t)-z_i) + \sum (\rho _z  \gamma _{\bar{y}}  + \rho _{\bar{z}} \bar{\gamma}_{\bar{y}}) z_i ^{'} (t)
\nonumber 
\end{eqnarray}
and similarly for the other quantities. Notice that $\vert z_i (t)-z_i \vert \leq E t^2$. For future reference let $\Vert \rho \Vert _{C^2} \leq \frac{\mathcal{P}}{10000}$ and $\max _i \vert z_i \vert = a$. Notice that by Cauchy's estimates $\vert z_i ^{'} \vert \leq 4(a+E)$ and $\vert z_i ^{''} \vert \leq 16(a+E)$ for $t<\frac{3}{4}$. For a system of IDE of the type
\begin{eqnarray}
w(s)-w &=& \displaystyle \int_{0}^s \mathbf{A}(t) w(t) dt \nonumber 
\end{eqnarray}
by the Gronwall inequality, $\vert w (s) \vert \leq \vert w\vert \exp(s\max_{[0,s]}\sqrt{n} \Vert \mathbf{A} \Vert)$. Without loss of generality, let $t<\min(3/4,\exp(-E\mathcal{P}/2)\frac{1}{\sqrt{100E\mathcal{P}}},\frac{}{})$ and $\mathcal{P}, E>1$. Therefore
\begin{eqnarray}
\Vert \gamma_y (t,1,y,\bar{y}) \Vert _{C^0} &\leq& \exp(\mathcal{P}Et^2) \leq \exp(\mathcal{P}E) \nonumber \\
\Vert \gamma_{\bar{y}} (t,1,y,\bar{y}) \Vert _{C^0} &\leq& \mathcal{P}Et^2\exp(\mathcal{P}Et^2) \leq \exp(\mathcal{P}E) \nonumber \\
\vert \gamma _y(s) - 1 \vert &\leq& 2\mathcal{P}Et^2 \exp(E\mathcal{P}) \leq \frac{1}{2} \nonumber \\
m &=& \frac{1}{2} \nonumber \\
\Vert \gamma _t (t,1,y,\bar{y} \Vert _{C^0} &<& \mathcal {P}E (a+E) \exp (\mathcal{P}E) \nonumber \\
\Vert \gamma _{tt} (t,1,y,\bar{y} \Vert _{C^0} &<& (\mathcal {P}E (a+E))^3 \exp (\mathcal{P}E)  \nonumber \\
\Vert \gamma_{yt} (t,1,y,\bar{y}) \Vert _{C^0} &<& (\mathcal{P}E(a+E))^2 \exp(\mathcal{P}E) \nonumber \\
\Vert \gamma_{\bar{y}t} (t,1,y,\bar{y}) \Vert _{C^0} &<& (\mathcal{P}E(a+E))^2 \exp(\mathcal{P}E)  \nonumber \\
\Vert \gamma_{ytt} (t,1,y,\bar{y}) \Vert _{C^0} &<& (\mathcal{P}E(a+E))^3 \exp(\mathcal{P}E) \nonumber \\
\Vert \gamma_{\bar{y}tt} (t,1,y,\bar{y}) \Vert _{C^0} &<& (\mathcal{P}E(a+E))^3 \exp(\mathcal{P}E)  \nonumber \\
C=2\frac{6M^4}{m^4}+1 &<& 200(\mathcal{P}E(a+E))^{12} \exp(4\mathcal{P}E) \nonumber 
\end{eqnarray}
Recall that $a = \max _i |z_i|$ (not $z_i(t)$ but $z_i(0)$), $E$ is that constant such that $\vert z_i (t) - z_i \vert \leq Et^2$, and $\mathcal{P} = 1000\Vert \rho \vert _{C^2}$. Note that $t$ has an additional condition in that $t<\min(\frac{3}{4},\exp(-E\mathcal{P}/2)\frac{1}{\sqrt{100E\mathcal{P}}})$ so as to ensure that indeed $\vert \gamma _y \vert \geq \frac{1}{2}$ and that $\vert v \vert < 1$. 

Until now the constants depended on $\max_{i} |z_{i}|$. However, note that the extremal map problem is invariant under M\"{o}bius transformations fixing the upper half plane. The constants $\delta$ and $C_r$ in the mapping theorem depend only on $r$ (where $r$ is the disk inside which these estimates are valid). For the polygon mapping problem, a-priori all the $z_{i}$s are on the real line, and three of them are $0,1$ and $\infty$. Assume that $\infty$ is the $n$th puncture (so $z_{n}= \infty$),and choose $z_{\text{min}} = \min_{i \neq n} z_{i}$ and $z_{\text{max}} = \max_{i\neq n} z_{i}$. We then find a M\"{o}bius transformation that maps $z_{\text{min}},z_{\text{max}}$ and $z_{n}$ to $0,1,\infty$, respectively. Now all the new punctures are in the interval $[0,1]$, and the mapping theorem provides absolute constants that do not depend on the punctures anymore.

\end{proofof}
}

We showed in \thmref{problemsrelation} how to reduce the polygon mapping problem to the punctured sphere problem. However, the above procedure can also be directly implemented on polygons, once we have the appropriate basis for the space of quadratic differentials. We give a simple well-known description of this basis in \subsecref{qdonpoly} in the appendix..

\section{Discretization of the procedure}
\seclabel{discrete}

Before we discretize the procedure, we give the properties of the mesh we work on as promised in \subsecref{ourresults}. Given an error tolerance $\delta$, let $\epsilon = O(\delta ^{6n-2}) $.

\begin{definition}\deflabel{canonicaltriangulation}[Canonical triangulation of size $\epsilon$] A canonical triangulation of size $\epsilon$, denoted as $\Delta_{\epsilon}$ is a set of vertices and edges $(V_{\epsilon},E_{\epsilon})$, with $z_{i} \in V_{\epsilon}$, satisfying the following.
\begin{CompactEnumerate}
\item It contains the edges and vertices of a regular polygon centered at $0$ and of diameter $O(\delta^{-1})$, and line segments joining the vertices of this polygon to $\infty$. 
\item Except for the line segments to $\infty$, all the other sides of the triangulation have Euclidean length at most $\epsilon$.
\item It contains the edges and vertices of regular polygons of $N = O(\delta^{(1-2n)})$ sides centered at the punctures of diameter $O(\delta)$, and lines joining the the vertices of these polygons to their centers, i.e. to the punctures. 
\item $\Delta_{\epsilon}$  is a refinement of $\Delta_{1/2}$.
\end{CompactEnumerate}
\end{definition}

We now describe what the algorithm \INEXT does, after which we prove \thmref{discrete1}.
\subsection{\INEXT and the proof of \thmref{discrete1}}
\seclabel{INEX}
We want to discretize the operator $\mathcal{P}(\mu)$ which returns $\nu$ with the least $L^{\infty}$ norm satisfying $\int_{R}\nu \phi_{i} = \int_{R} \mu \phi_{i}$ for all $\phi_{i}$ in Equation~\eqreff{puncturediff}. Note that the starting $\mu$ is piecewise constant at the start of every iteration.
\begin{observation}
The integral of $\phi_{i}$ over any triangle $t_{j}$ can be computed analytically.
\end{observation}
We provide this formula (that involves taking the logarithm of a complex number) in \subsecref{formula_integral}. 

We approximate $\nu$ by piecewise constant Beltrami differentials. One can easily see then that the constraints for infinitesimal equivalence become linear constraints of the form $Ax=b$, where $A$ is the matrix whose $(i,j)$th entry equals $\int_{t_{j}} \phi_{i}$, $x$ is the vector of unknown values of the piecewise constant $\nu$ on a triangle, and $b$ is the vector of $\int_{t_{j}} \mu_{j} \phi_{i}$, where $\mu_{j}$ is the value of $\mu$ on triangle $t_{j}$.

If $A$, $x$ and $b$ are real, an $L^{\infty}$ minimization can be formulated as a linear program. In our case, we break the vectors and matrices into their real and complex parts, and then we can formulate the program as a quadratically constrained quadratic program. Although in general they are NP-hard to solve, we show that our program involves positive semi-definite matrices; and it is known that such instances can be solved in polynomial time using interior-point methods \cite{pardalos2002handbook}. Details are in \subsecref{inextdetails}. 

%Using the observation that the true $\nu$ is of the form $k \bar{\phi}/|\phi|$, we formulate the minimization program as a quadratically constrained quadratic program (QCQP), which can be solved in polynomial (in the bit complexity of the input) time by interior point methods. We give the details in appendix \subsecref{inextdetails}. Thus
\begin{lemma}\lemlabel{def_inext}[\INEXT]
There exists an algorithm \INEXT that, given a piecewise constant $\mu$ on $\Delta_{\epsilon}$ returns a piecewise constant $\hat{\nu}$ such that $\max_{t_{j}} \hat{\nu}(t_{j}) \leq  \max_{t_{j}} \beta(t_{j})$, where $\beta$ is any piecewise constant (on $\Delta_{\epsilon}$) Beltrami differential that is infinitesimally equivalent to $\mu$.
\end{lemma} 

With this, we are now in a position to prove \thmref{discrete1}.

\textbf{Proof sketch for \thmref{discrete1}:}
The main idea is to use $\nu = \mathcal{P}(\mu)$ and produce a piecewise constant $\nu_{p}$, which does not satisfy the integral constraints, but the error involved can be estimated. We then add a piecewise constant differential to $\nu_{p}$ to produce $\tilde{\nu}_{p}$ that is  in the same infinitesimal class as $\mu$ and whose dilatation is close to that of $\nu$. This proves that $\hat{\nu}$ (whose dilatation is smaller than that of $\tilde{\nu}_{p}$) satisfies the desired requirement.

\textlater{\begin{proofof}{\thmref{discrete1}}
In whatever follows we assume that $\vert z_i \vert \leq A$. Let $\tilde{\tilde{k}}=\min(\vert z_i -z_j \vert/200, \vert z_i - 1 \vert /200, \vert z _i \vert/200) $. Here we prove that $\Vert \nu \Vert - \Vert \nu_{*} \Vert < \delta$. The strategy is to produce a piece-wise constant Beltrami coefficient $\tilde{\nu}_{p}$ whose norm is $\delta$ close to $\nu_{*}$ and satisfies the integral constraint exactly. To do this we first produce a piece-wise constant Beltrami $\nu_{p}$ whose norm is equal to that of the infinitesimally extremal one $\Vert \nu_{*}\Vert$ and $\vert \int (\nu_{p} - \nu_{*}) \phi_k \vert < \delta \ \forall \ k$. We claim that we may produce $\tilde{\nu}_p$ by adding terms of magnitude at most $\delta$ to $\nu_{p}$ so as to make sure that the integral constraint is obeyed exactly. Indeed, by the hypothesis on the canonical triangulation we see that one may choose sufficiently many triangles of size $\frac{1}{2}$ and solve linear equations that determine the constants to be added on these triangles so as to satisfy the integral constraint. Thus the problem is reduced to finding $\nu_p$ and proving the error estimate. 

The infinitesimally extremal $\nu_{*} = k_{*} \frac{\vert \phi \vert}{\phi} = k_{*}\frac{\vert \sum c_ k \phi _k \vert}{\sum c_k \phi_k } = k_{*}\frac{\vert \sum \tilde{c}_ k \phi _k \vert}{\sum \tilde{c}_k \phi_k } $ where $\vert \tilde{c}_k \vert \leq 1$ with $\vert \tilde{c}_1 \vert =1$. Then 
\begin{gather}
\nu_{*} = k_{*} \frac{\vert\sum \tilde{c}_k (z-z_1)\ldots (z-z_{k-1})(z-z_{k+1})\ldots \vert}{\sum \tilde{c}_k (z-z_1)\ldots (z-z_{k-1})(z-z_{k+1})\ldots}\frac{z(z-1)(z-z_1)\ldots}{\vert z(z-1)(z-z_1)\ldots \vert} \nonumber \\
= k_{*}  \frac{\vert\sum \tilde{c}_k h_k(z)\vert}{\sum \tilde{c}_k h_k(z)\ldots}g(z) \nonumber
\end{gather}

We first approximate $h_k$ and $g$ by piece-wise constant functions with an error of at most $\beta$ within the large polygon. Replace the values of $h_k$, $g$ by their values at the centers of the triangles in the triangulation. This gives us $\nu_{p}$ in the large polygon. Outside it we set $\nu_p$ to $0$.

The error thus caused for $h_k$ is less than $\epsilon 2^n R^n$. For $g$ the error is more subtle. In the small polygons around the punctures, if $N > \frac{12\pi}{\beta}$ then the error caused in $g(z)$ is less than $\frac{\beta}{3}$ in those regions. For future use it is useful to note that $\frac{y_2}{\vert y_2} - \frac{y_1}{\vert y_1} \leq \vert y_2 - y_1 \vert  \int _{0} ^1 \vert \frac{1}{ty_2 + (1-t) y_1} \vert dt$ and the same inequality holds for $\frac{\vert y \vert}{y}$. Outside the polygons the error in $g$ is less than $\frac{\sqrt{\epsilon} 2^n R^n}{1-\sqrt{\epsilon}2^n R^n}$. Since the radius of the polygons is $\tilde{\kappa}$, the error in $\frac{\vert \sum \tilde{c}_k h_k \vert}{\sum \tilde{c}_k h_k}$ outside the polygons is less than $\frac{\epsilon 2^n R^n}{(\tilde{k})^{n-1} - \epsilon 2^n R^n}< \beta$. Notice that $\epsilon < \min(\tilde{k} ^{2n},\frac{\beta ^2}{2^{2n}R^{2n}})$. 
Next we want to estimate $I_k = \int \displaystyle  \frac{\vert \eta_k (\nu _{*}-\nu_{p} ) \vert}{\vert z(z-1)(z-z_k ) \vert} = \int U_k $. Indeed,
\begin{gather}
I_k \leq 4k_{*} \int _{\vert z \vert \geq 2R} \frac{1}{r^2} dr d\theta + \int _{\vert z \vert \leq 2R} U_k \nonumber \\
\leq \frac{k_{*}\delta}{2}  + J_k \nonumber
\end{gather}
where $2\vert z \vert ^3 \geq \vert z(z-1) (z-z_k) \vert \geq \frac{\vert z \vert^3}{2}$ if $R > \max{\frac{4}{\delta},1,2^n (1+A)^n}$. Now we estimate $J_k$. 
\begin{gather}
J_k = \displaystyle \int _{\vert z \vert \leq 2R \cap \vert \sum \tilde{c}_k h_k \vert > \kappa} U_k + \int _{\vert z \vert \leq 2R \cap \vert \sum \tilde{c}_k h_k \vert > \kappa} U_k \nonumber \\
< k_{*} \frac{\beta}{\kappa - \beta}  + k_{*}\beta +  \displaystyle \sum _{\xi}\int _{\vert z \vert \leq 2R \cap \vert \sum \tilde{c}_k h_k \vert > \kappa \cap \vert z - z_{\xi} \vert < \tilde{\kappa}} U_k \nonumber + \int _{\vert z \vert \leq 2R \cap \vert \sum \tilde{c}_k h_k \vert > \kappa \cap \vert z - z_{\xi} \vert > \tilde{\kappa} \ \forall \ \xi}\\
< k_{*} \frac{\beta}{\kappa - \beta}  + k_{*}\beta +  k_{*} (1+A)^2 \frac{12 \pi \tilde{k}}{(\tilde{\tilde{\kappa}})^2} + k_{*}(A+1)^2 \frac{(\pi \kappa ^{2/(n-1)})^n}{\tilde{\kappa} ^3} 
\end{gather}
where the last estimate is obtained by remembering that $\sum \tilde{c}_k h_k = (z-\lambda _1) \ldots$ for some $\lambda _i$ within the big polygon. One may choose $\tilde{\kappa}$, $\kappa$ and $\beta$ so as to estimate $I_k$ by $\delta$. Indeed,
\begin{gather}
\tilde{\kappa} < \frac{\delta (\tilde{\tilde{\kappa}})^2}{1000 (1+A)^2} \nonumber \\
\kappa < \left ( \frac{\delta ^4 (\tilde{\tilde{\kappa}})^6}{10^{10} (1+A)^8} \right ) ^{(n-1)/2} \nonumber \\
\beta < \delta \frac{\kappa}{6+\delta} \nonumber \\
< O(\delta ^{2n-1}) \nonumber
\end{gather}

\end{proofof}
}

\subsection{\EXTREMAL and the proof of \thmref{discrete2}}

Apart from the subroutines \INEXT, \PIECEWISECOMP and \PIECEWISEINV, we will require three more subroutines to discretize our procedure.

\begin{definition}\deflabel{triang}[Subroutine:\TRIANG]

\noindent\textbf{Input:} a set of points $\mathcal{S}$, a size $M$, and a triangulation $\Delta_{\epsilon}$.

\noindent\textbf{Output}: A triangulation $\Delta_{\epsilon^{'}}$ of the given size $M$ containing $\mathcal{S}$ such that $\Delta_{\epsilon^{'}}$ is a refinement of $\Delta_{\epsilon}$.
\end{definition}

\begin{definition}\deflabel{beltram}[Subroutine: \BELTRAMI]

\noindent\textbf{Input:} A triangulation $\Delta_{\epsilon}$ of the plane, a piecewise constant Beltrami coefficient $\mu$, and error tolerance $\delta$.

\noindent\textbf{Output:} A triangulation $\Delta_{\epsilon}^{'}$ that is a refinement of $\Delta_{\epsilon}$, and the images $\hat{f}(v_{i})$ of the vertices $v_{i} \in \Delta_{\epsilon}^{'}$ such that $|f^{\mu}(v_{i})-\hat{f}(v_{i})|<\delta$.
\end{definition}

\begin{definition}\deflabel{vecfield}[Subroutine: \VECTFIELD]

\noindent\textbf{Input:} A $C^k$ ($k$ sufficiently large, e.g. $k>10$) vector field $X$ (written as a formula in terms of elementary functions), a triangulation $\Delta_{\epsilon}$, and an error tolerance $\delta$.

\noindent\textbf{Output:}A triangulation $\Delta_{\epsilon}^{'}$ that is a refinement of $\Delta_{\epsilon}$, the images of $v_{i} \in \Delta_{\epsilon}$ up to error $\delta$ under a $C^k$ diffeomorphism $\gamma_{x}$ corresponding to the flow along $X$, and  a piecewise smooth Beltrami coefficient that approximates the  one up to error $\delta$.

\end{definition}

\paragraph{Implementing \TRIANG, \BELTRAMI and \VECTFIELD:} We outline ways to implement the above three subroutines:
\begin{CompactEnumerate}
\item Given a set of $n$ points, we can obtain the Delaunay triangulation in $O(n \log n)$ time. While implementing \TRIANG, we first compute the Delaunay triangulation of all the points falling inside a triangle of the given triangulation. The we connect the vertices on the convex hull of such a set of points to one of the three vertices of the triangle they lie in. If this complete triangulation is not yet size $M$, we make the mesh denser by adding points as in \cite{Ruppert95Delaunay} (points are added to either the circumcenters of triangles or mid-points of edges), until we reach the desired size.

\item The solution to the Beltrami equation for $\mu$ can be expressed as a series of singular operators applied to $\mu$. There are many efficient algorithms and implementations (\cite{Daripa1993355},\cite{doi:10.1137/050640710}) existing for \BELTRAMI. Most of them can bound the $l^{p}$ norm of the error, but the methods in \cite{Daripa1993355} can be used to bound the $L^{\infty}$ error too \cite{daripa}. 

\item The idea of deforming a surface by a vector field has been applied extensively in computer graphics. We refer the reader to \cite{lddmm} for an implementation. 
\end{CompactEnumerate}

\paragraph{Description of \EXTREMAL:} The algorithm summarized below is based on \subsecref{summarycontinuous}.
\begin{CompactItemize}
\item Use \TRIANG to produce a triangulation of size required by \INEXT to run within an error of $\delta^{10}$.
\item Loop $i=1$ to $N$ where $N$ is the number of iterations in \thmref{continuous2} to produce the result within an error of $\delta/2$.
\begin{CompactEnumerate}
\item Use \INEXT to produce $\nu_i$ from $\mu_{i}$ within an error of $\delta ^{10}$. If $\nu _{i} = \mu_{i}$ then stop. 
\item Find $t_i$ by Equation~\eqreff{choosing_t}, using $k_0$ as $\Vert \nu _i \Vert _{\infty}$. 
\item  Invoke \BELTRAMI for the coefficient $t_{i}(\mu_{i}-\nu_i)$ to find the images of the punctures within an accuracy of $t_i ^3$. 
\item Define the vector field $X$ as in the continuous construction using a piecewise polynomial version of the bump function (that is $C^{10}$ for instance). Then call \VECTFIELD to find a piecewise constant Beltrami coefficient up to an error of $t_i^3$. 
\item Use \PIECEWISECOMP to compose the Beltrami coefficients of step $3$ and step $4$ within an error $(\Vert \mu_{i} \Vert - \Vert \nu_{i} \Vert)^5$ for the Beltrami coefficient and $\delta/i^{2}$ for the q.c.h. 
\item Use \PIECEWISEINV to find the Beltrami coefficient of the inverse of the q.c.h of step $5$, up to the same error as that in step $5$.
\item Call \PIECEWISECOMP to compose $\mu_{i}$ and the Beltrami coefficient of step $6$ to form $\mu_{i+1}$ (up to the same error as that in step $5$).
\end{CompactEnumerate}
\end{CompactItemize}

The algorithm terminates by producing $\mu_{N}$. The proof of \thmref{discrete2} is similar to that of \thmref{continuous2}.

\section{Discussions}
\seclabel{discussion}

We conjecture our algorithm to run in polynomial time. This is evidenced by the fact that 1) the number of iterations is a polynomial in $1/\epsilon$, 2) \INEXT (quadratic program), \TRIANG, \BELTRAMI and \VECTFIELD run in polynomial time, and 3) we expect the existence of polynomial time subroutines \PIECEWISECOMP and \PIECEWISEINV.

Open problems abound. Apart from improving (complexity and approximation) the algorithm we propose, the extremal map problem can be further explored in many directions.
\begin{CompactEnumerate}
\item Most of the ideas presented here (notably \lemref{reductionlemma}) can be used to envision an algorithm for computing Teichm\"{u}ller maps between arbitrary (finite analytic type) Riemann surfaces. The problem is challenging for multiple reasons---for instance, an explicit basis of holomorphic quadratic differentials may not be available. 

\item The authors feel that building a discrete version of Teichm\"{u}ller theory would be an important achievement. Given a triangulated Riemann surface, defining a discrete analog of dilatation that gives nice results (e.g. existence and uniqueness) about the extremal map would be the next step in this direction.

\item Most of the surfaces we see in everyday life can be regarded as Riemann surfaces. Being able to compute the "best" angle-preserving map between them is certainly of theoretical and practical interest. Our current efforts are aimed at being able to visualize geodesics in Teichm\"{u}ller space. Seeing the base polygon (or Riemann surface) morphing (similar to what was accomplished in \cite{Sharon:2006:AUC}) into the target polygon (surface) under the solution to $t \mu_{*}$, would give us an idea of how shapes actually change while following a geodesic in this moduli space.
\end{CompactEnumerate}

\noindent\textbf{Acknowledgements:} The first author would like to thank Frederick Gardiner, Christopher Bishop, Irwin Kra and Joe Mitchell for numerous discussions and helpful suggestions.
\newpage

\bibliographystyle{plain}

\newpage
\begin{center}
{\Large{\textbf{Appendices}}}
\end{center}

\section{Appendix $1$: Quasiconformal maps and essential theorems from Teichm\"{u}ller Theory}
\seclabel{app1}

All of the material in this section is classical and can be found in books on complex analysis and Riemann surfaces, such as  \cite{Ahlfors_2006,Kra04,Gardiner,forster}.

\subsection{Riemann mapping theorem and Riemann surfaces}
\subseclabel{rmtrs}

\begin{theorem}[Riemann Mapping]
\thmlabel{riemannmapping}
Let  $\Omega$  be a simply connected domain in the complex plane $\mathbb{C}$, not equal to the entire complex plane. Then there exists a biholomorphic map $f : \mathbb{D} \longrightarrow \Omega$. Further, $f$ is unique up to composition by a M\"{o}bius transformation.
\end{theorem}

That $f$ is biholomorphic implies it is conformal. One can therefore state as a corollary that any two simply connected domains in $\mathbb{C}$ (not equal to $\mathbb{C}$) can be mapped conformally and bijectively to each other.

\paragraph{Riemann surface}

Let $M$ be a two dimensional real manifold. A complex chart on $M$ is a homeomorphism $\phi$ from an open subset $A \subset R$ to an open subset $B \subset \mathbb{C}$. Let $\phi_{1}:A_{1} \rightarrow B_{1}$ and  $\phi_{2}:A_{2} \rightarrow B_{2}$ be two complex charts. $\phi_{1}$ and $\phi_{2}$ are said to be compatible if the map
\[ \phi_{2} \circ {\phi_{1}}^{-1} : \phi_{1}(A_{1} \cap A_{2}) \rightarrow \phi_{2}(A_{1} \cap A_{2})\] is biholomorphic.

A complex atlas on $M$ is a system of charts which cover $M$, and in which any two charts are compatible. Two complex atlases are regarded equivalent if all charts in the union of the atlases are pairwise compatible.

\begin{definition}[Riemann surface]
\deflabel{riemannsurfacedef}
A Riemann surface $R$ is a pair $(M, \sigma)$, where $M$ is a connected two-manifold and $\sigma$ is an equivalence class of complex atlases on $M$.
\end{definition}

Examples of Riemann surface include the complex plane, domains in the complex plane, the Riemann sphere $\hat{\mathbb{C}}$ and all Riemannian manifolds (oriented two-manifolds with a Riemannian metric).

Given two Riemann surfaces $M$ and $N$, a map $f:M\to N$ is \emph{conformal} if its restriction on any local conformal parameters is \emph{holomorphic}. Geometrically, a conformal map preserves angles, and transforms infinitesimal circles to infinitesimal circles, as shown in Figure~\ref{fig:qcillustration2} frame (a),(b) and (c).

\subsection{Quasiconformal maps}
\subseclabel{qcmapsapp}

\begin{figure*}[h]
\centering
\begin{tabular}{ccccc}
\includegraphics[width=0.15\textwidth]{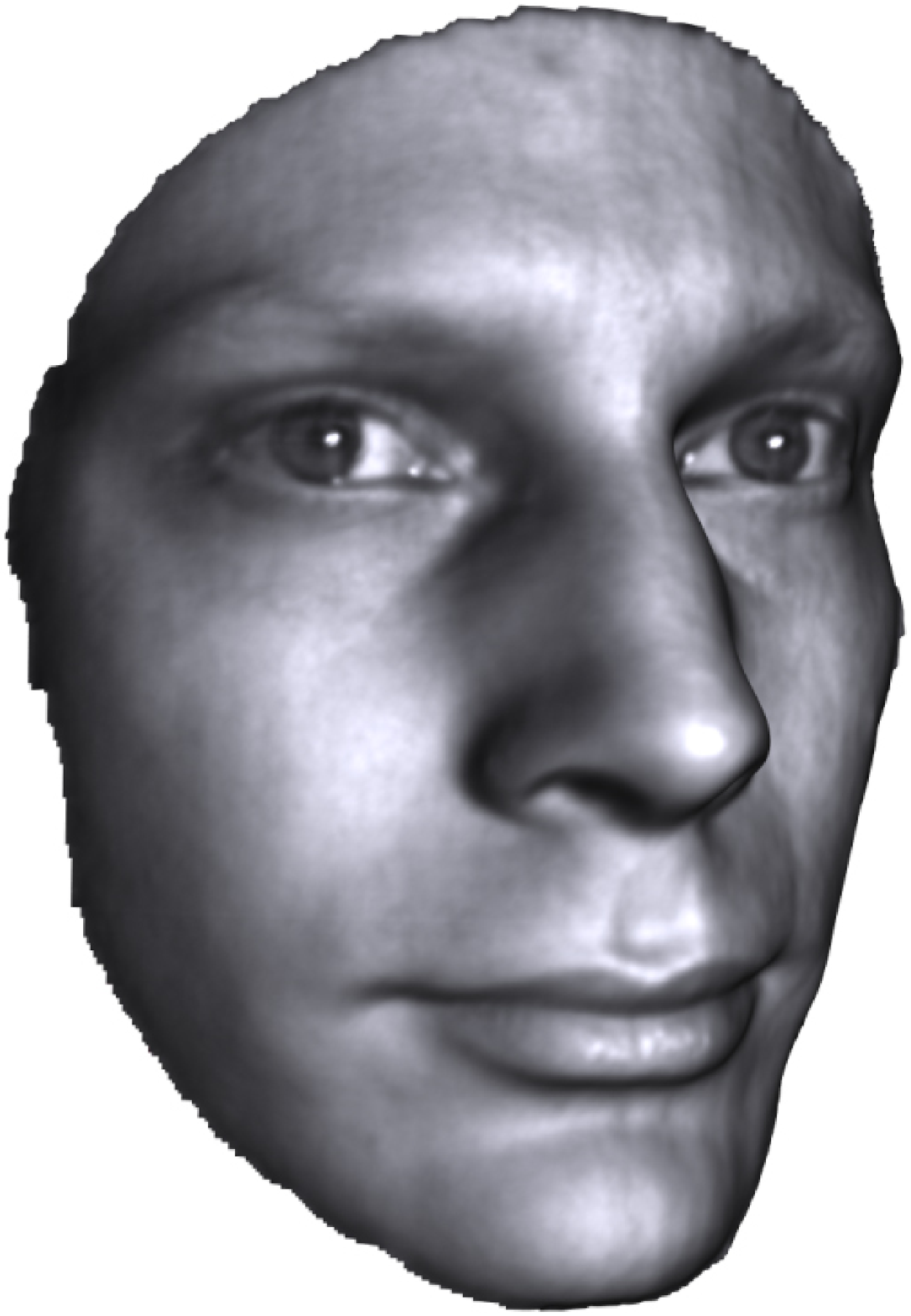}&
\includegraphics[width=0.22\textwidth]{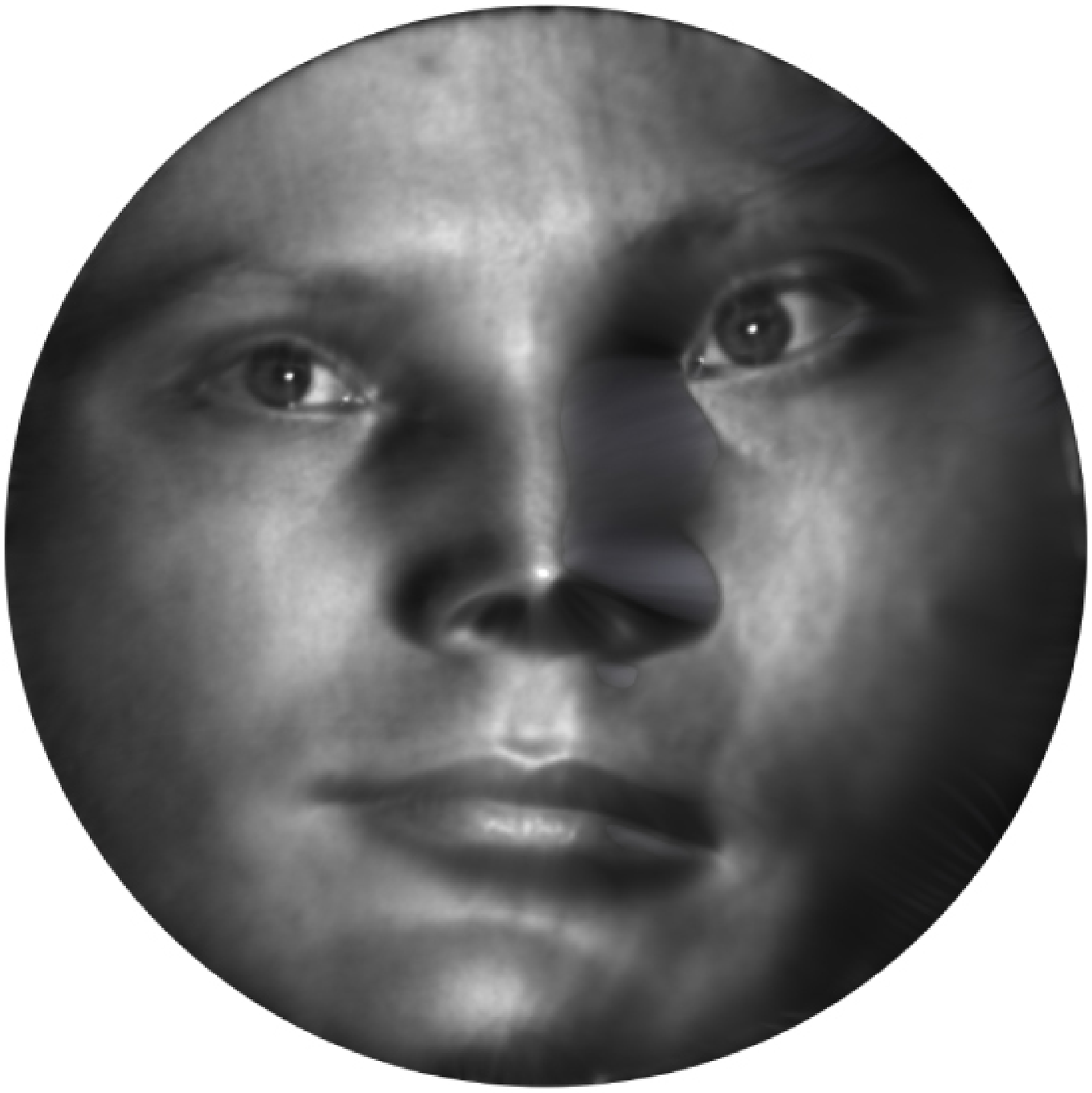}&
\includegraphics[width=0.15\textwidth]{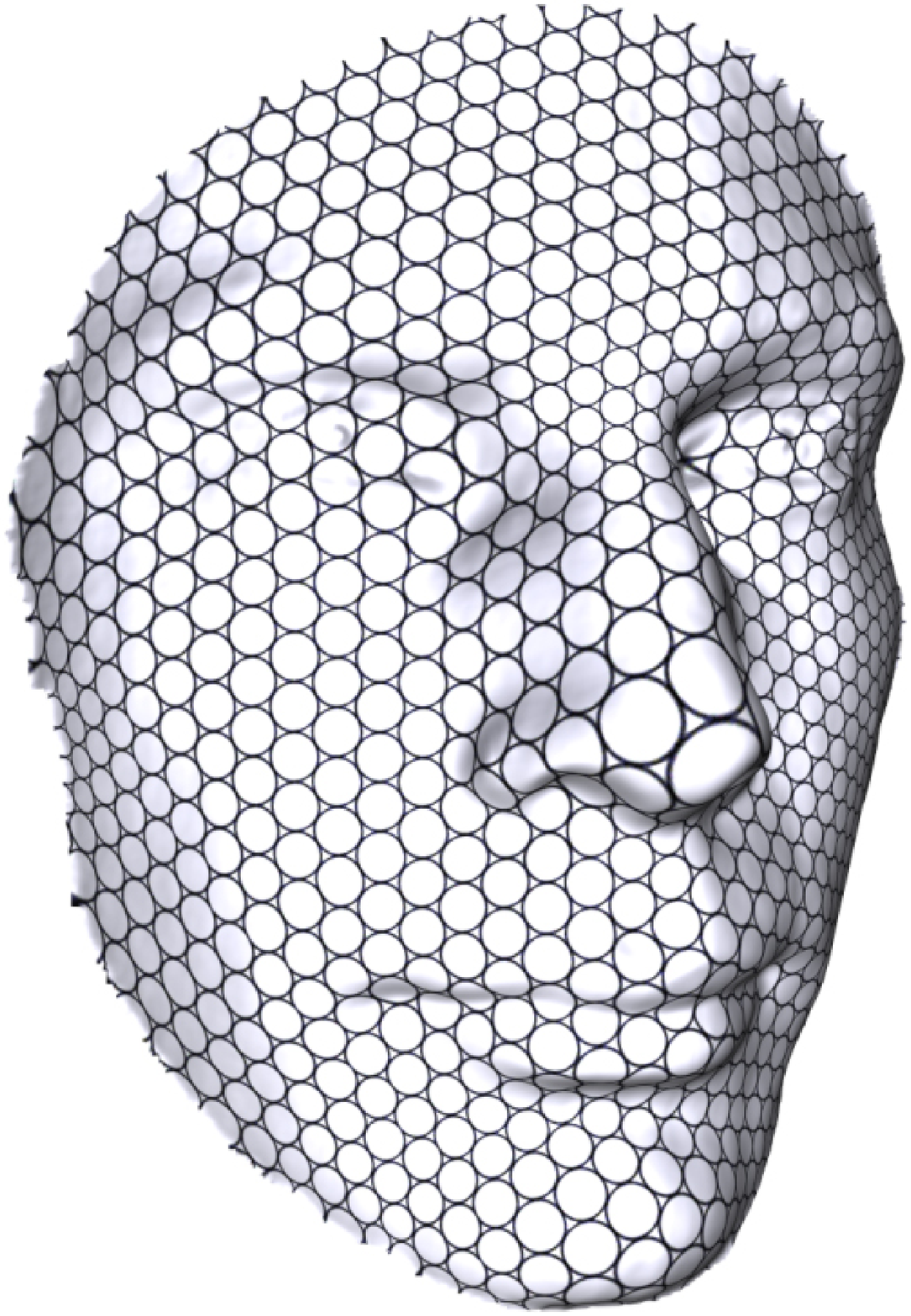}&
\includegraphics[width=0.22\textwidth]{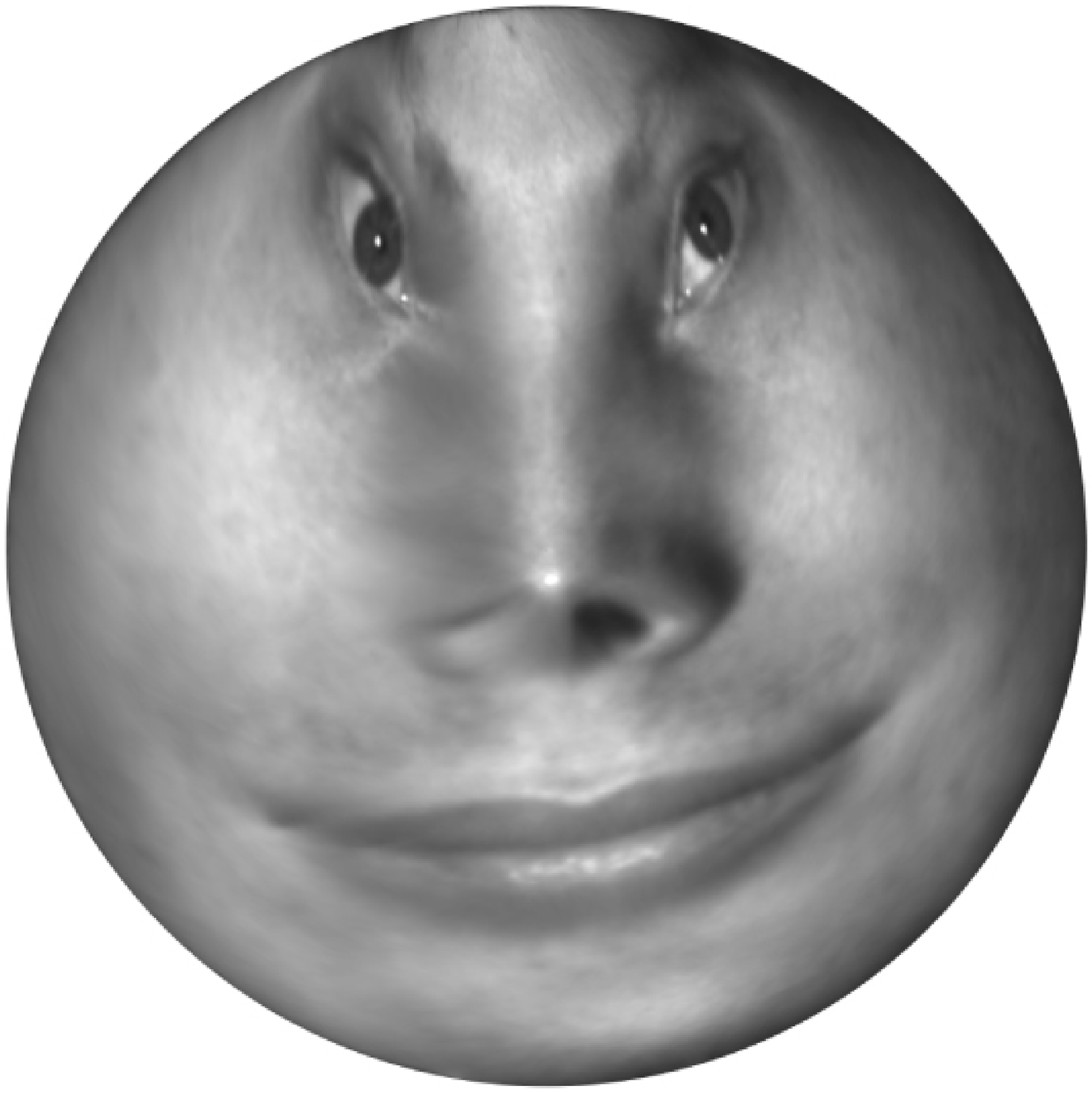}&
\includegraphics[width=0.15\textwidth]{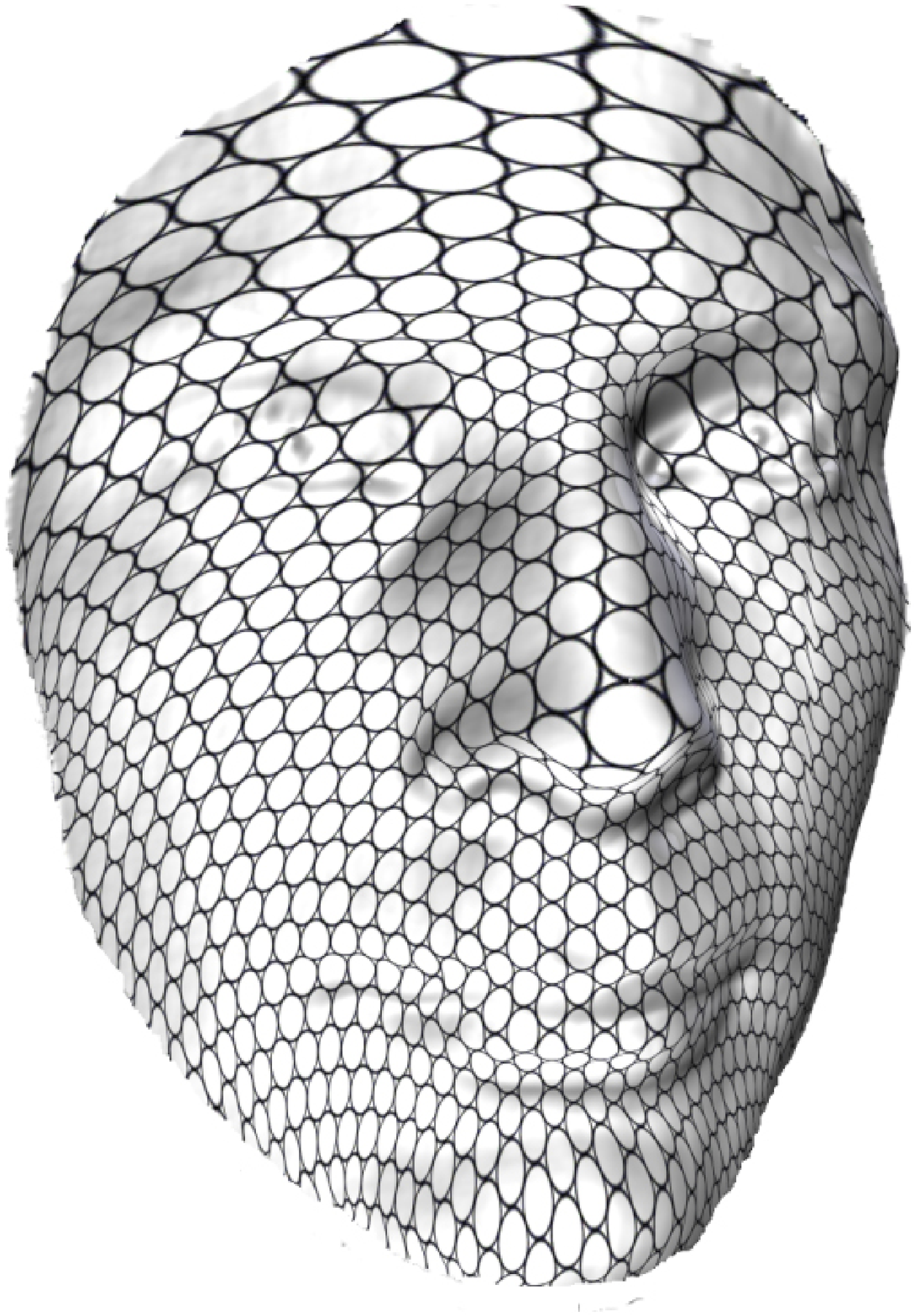}\\
(a)Original surface&(b) Conformal & (c) Conformal &(d) Qc mapping  & (e) Qc mapping\\
&mapping & mapping & & \\
\end{tabular}
\caption{Conformal and quasiconformal mappings from a human face surface to the planar disk. \label{fig:qcillustration2}}
\end{figure*}

\begin{figure*}[h!]
\centering
\begin{tabular}{c}
\includegraphics[width=0.35\textwidth]{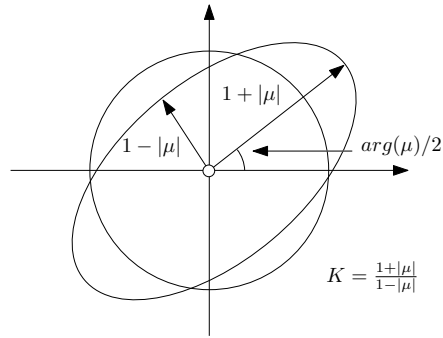}\\
\end{tabular}
\caption{Beltrami coefficient. \label{fig:BeltramiCoefficient}}
\end{figure*}

A generalization of conformal maps are \emph{quasiconformal} maps, which are orientation preserving homeomorphisms between Riemann surfaces with bounded conformality distortion, in the sense that their first order approximations takes small circles to small ellipses of bounded eccentricity, as shown in Fig.\ref{fig:qcillustration2} frame (d) and (e). Mathematically, $f \colon \mathbb{C} \to \mathbb{C}$ is quasiconformal provided that it satisfies the Beltrami equation:
\begin{equation}\label{beltramieqt}
f_{\bar{z}} = \mu(z) f_z.
\end{equation}
\noindent for some complex-valued function $\mu$ satisfying $||\mu||_{\infty}< 1$. $\mu$ is called the \emph{Beltrami coefficient}, and is a measure of the non-conformality of $f$.  In particular, the map $f$ is conformal around a small neighborhood of $p$ when $\mu(p) = 0$.
As shown in Figure~\ref{fig:BeltramiCoefficient}, the orientation of the ellipse is double the argument of $\mu$. The \emph{dilatation} of $f$ is defined as the ratio between the major axis and the minor axis of the infinitesimal ellipse. The maximal dilatation of $f$ is given either by:

\begin{equation}
k_{f} =||\mu_{f}||_{\infty},
\end{equation}
or by
\begin{equation}
K(f) = \frac{1+||\mu||_{\infty}}{1-||\mu||_{\infty}}.
\end{equation}

A homeomorphism with dilatation less than or equal to $K$ is called a \emph{K-quasiconformal mapping}.

\subsection{Quadratic differentials}
\label{quadraticdifferentials}

\begin{figure*}[h!]
\centering
\begin{tabular}{cccc}
\includegraphics[width=0.25\textwidth]{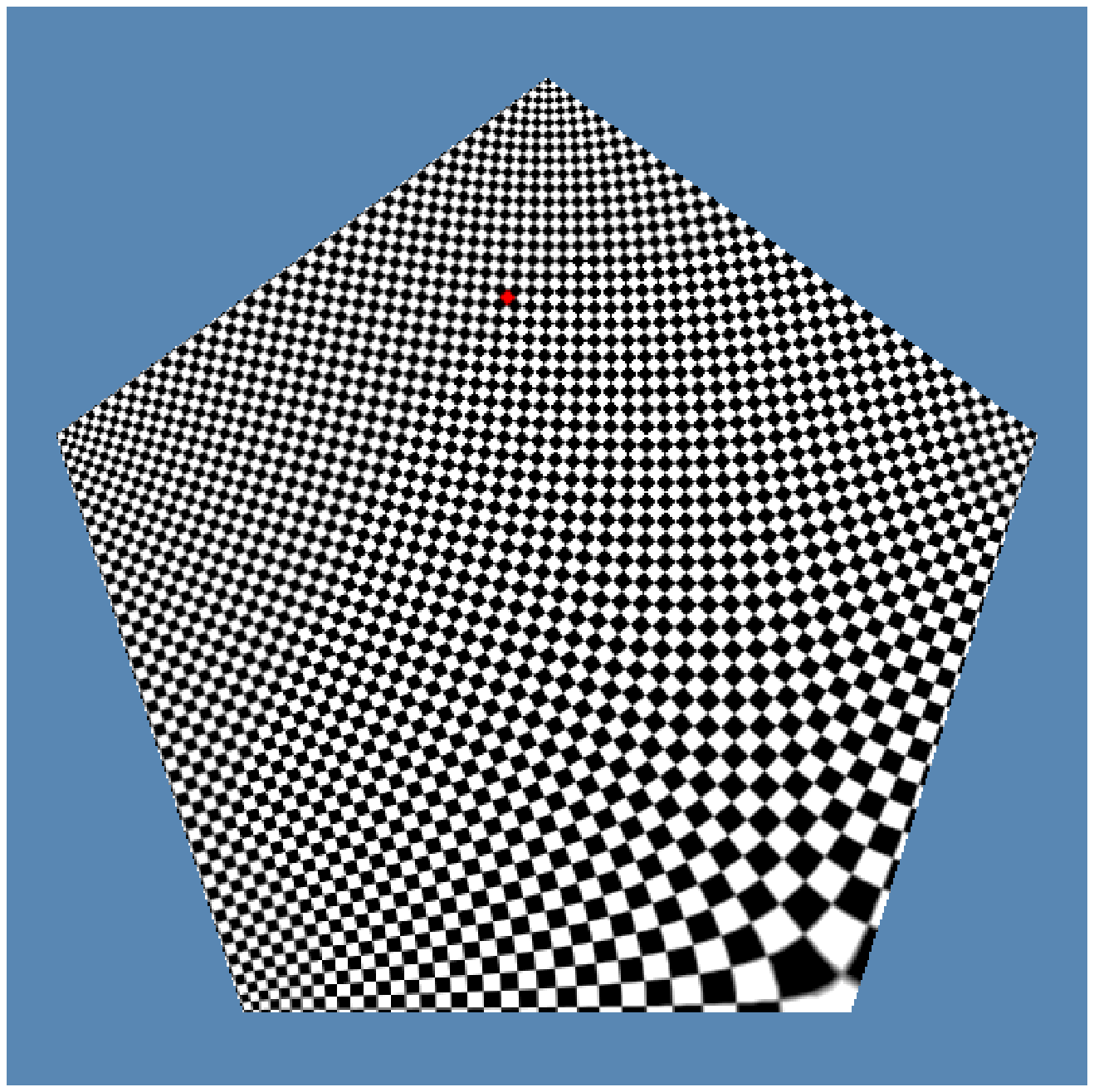}&
\includegraphics[width=0.25\textwidth]{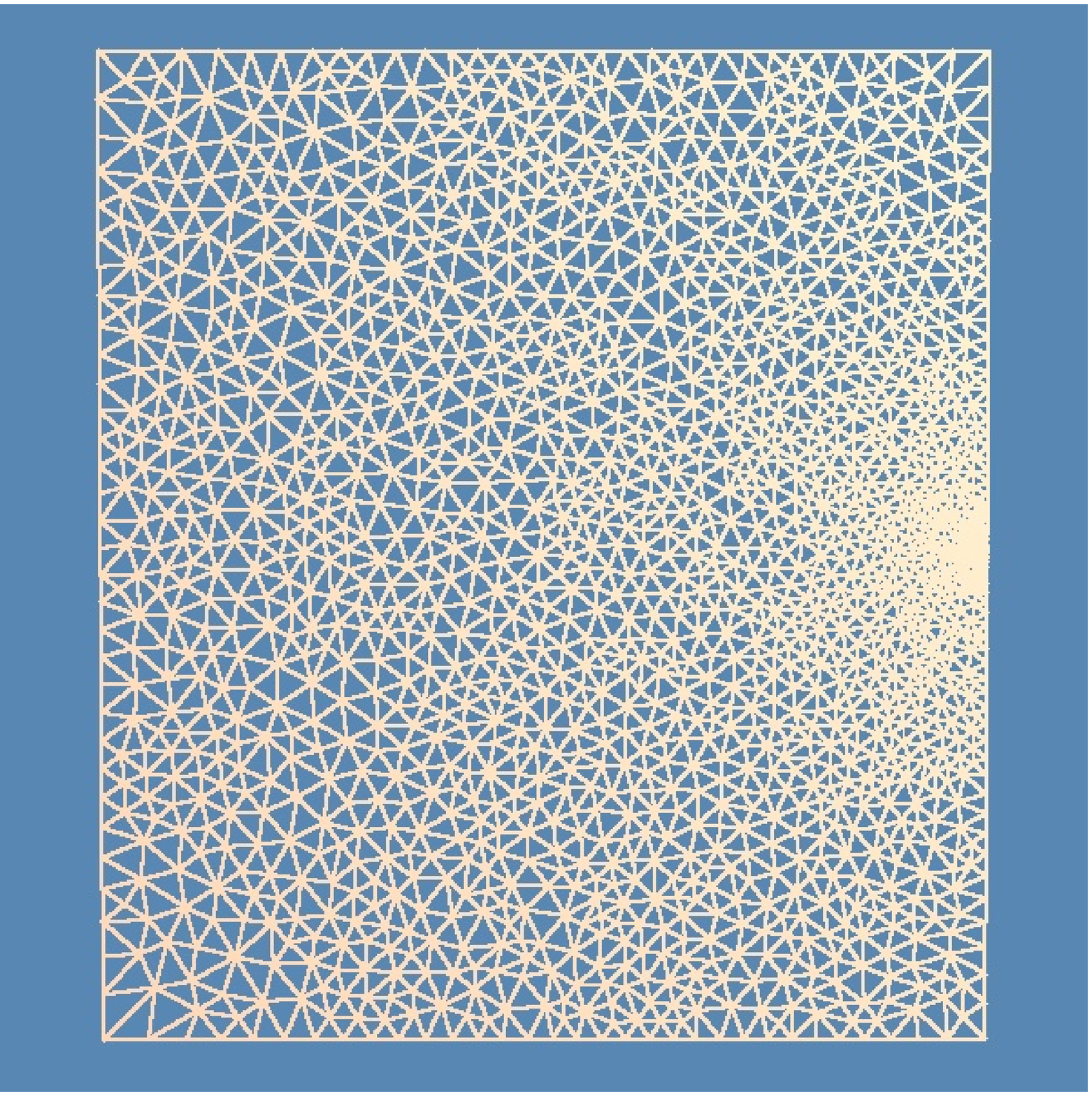}&
\includegraphics[width=0.25\textwidth]{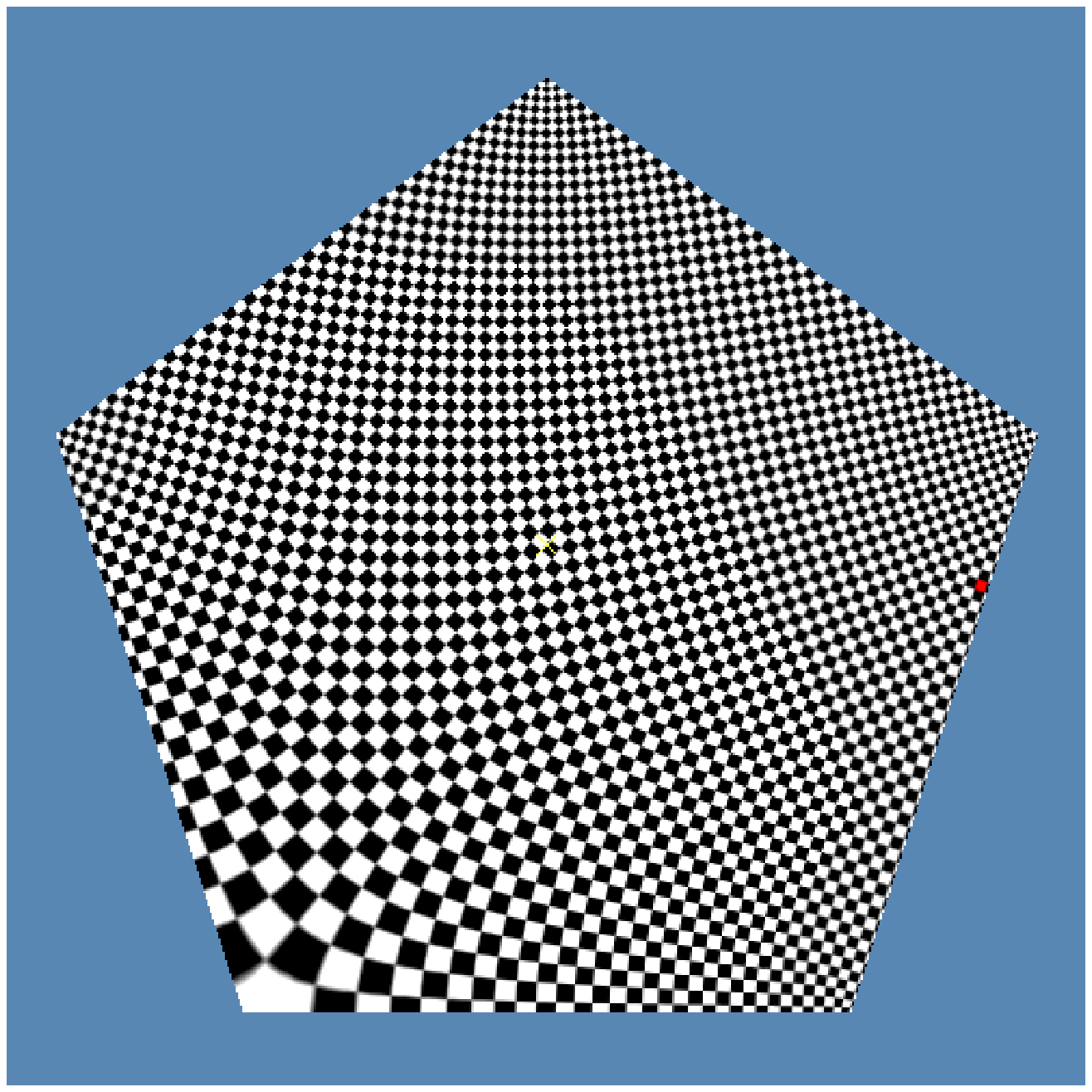}&
\includegraphics[width=0.25\textwidth]{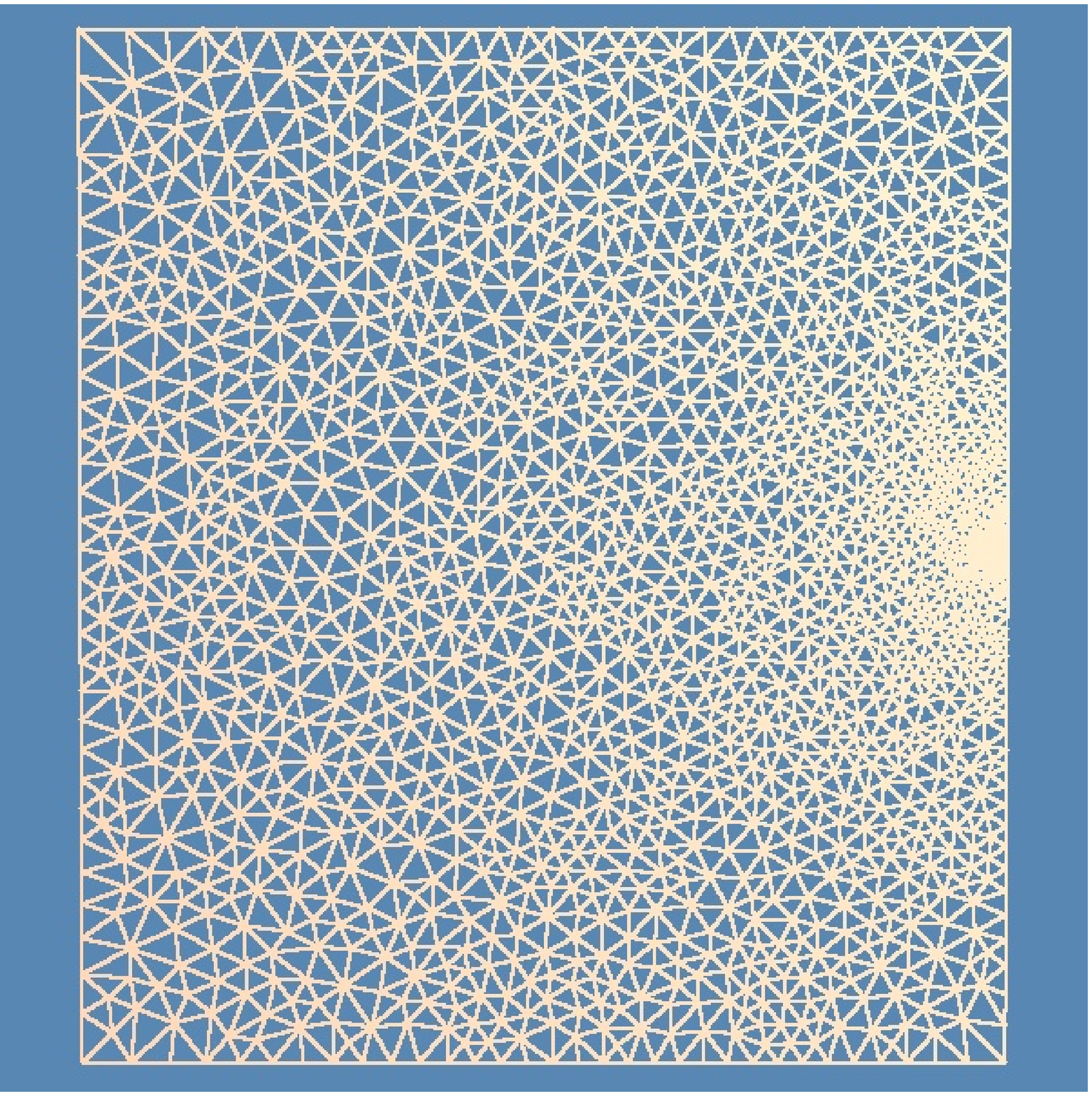}\\
(a)$\phi_1$ &(b)$R_1$ & (c) $\phi_2$ & (d) $R_2$\\
\end{tabular}
\caption{Holomorphic quadratic differential bases on a pentagon. \label{fig:pentagon_quad_diff_base}}
\end{figure*}

\begin{figure*}[h!]
\centering
\begin{tabular}{cccc}
\includegraphics[width=0.25\textwidth]{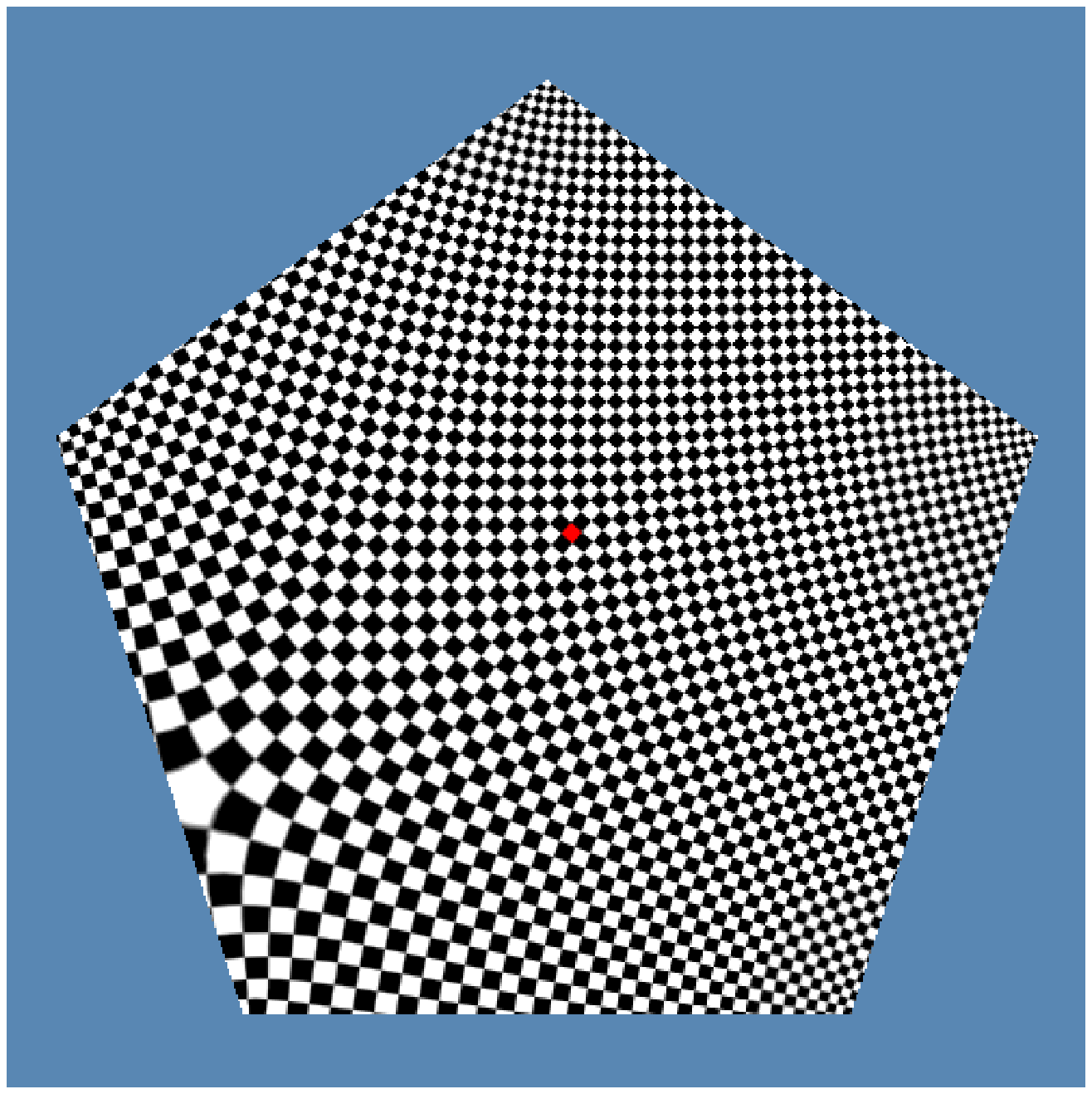}&
\includegraphics[width=0.25\textwidth]{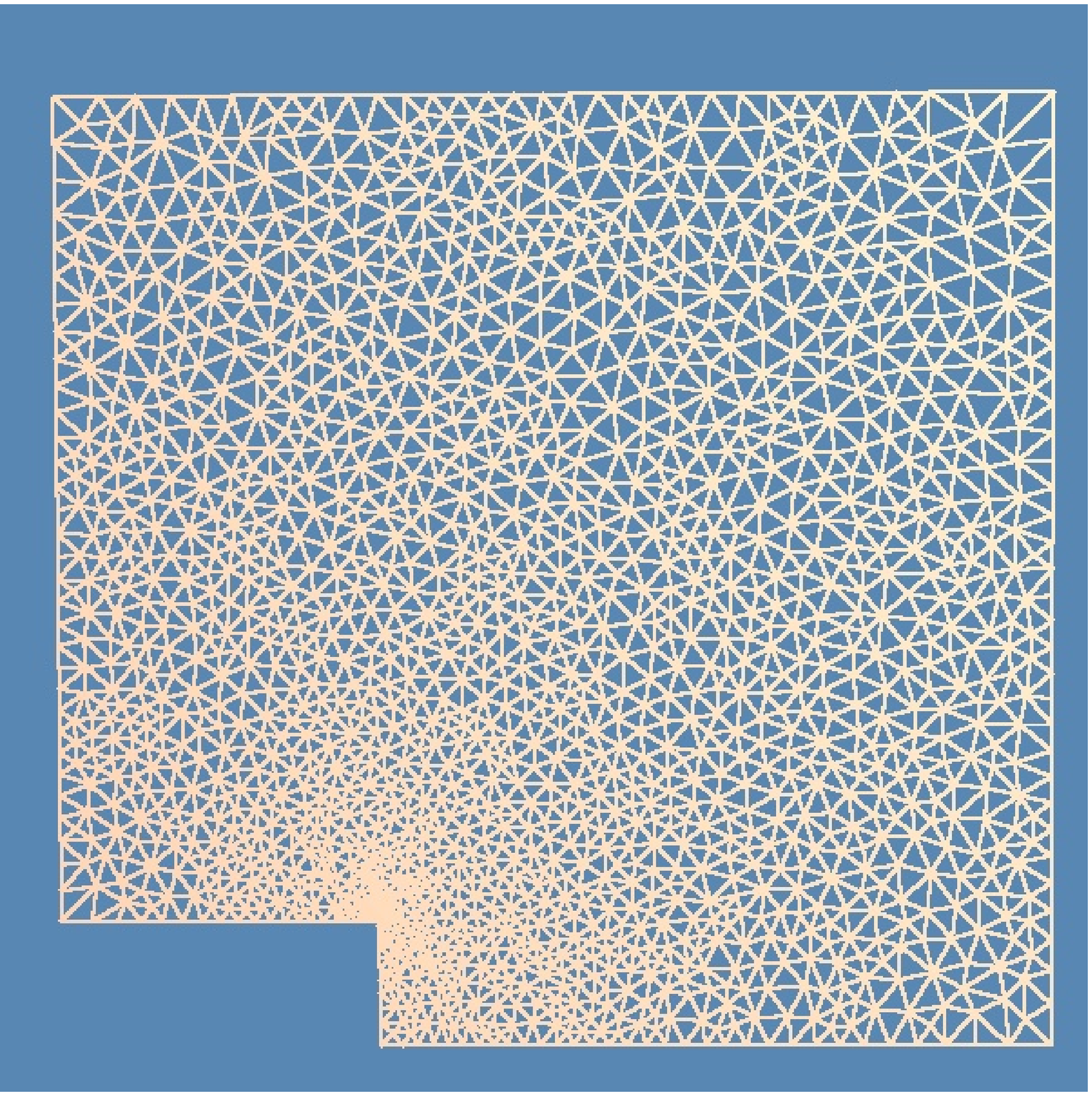}&
\includegraphics[width=0.25\textwidth]{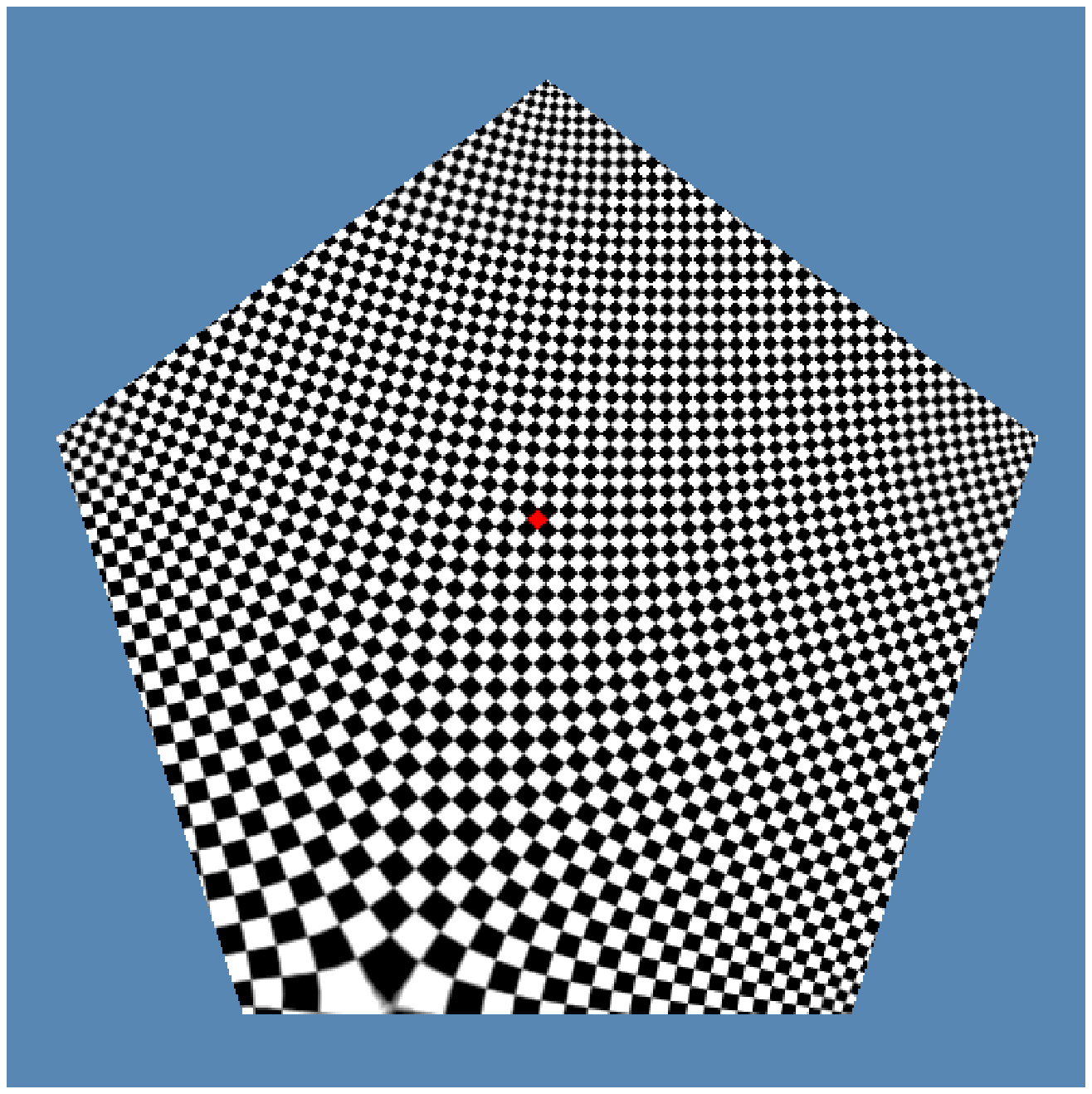}&
\includegraphics[width=0.25\textwidth]{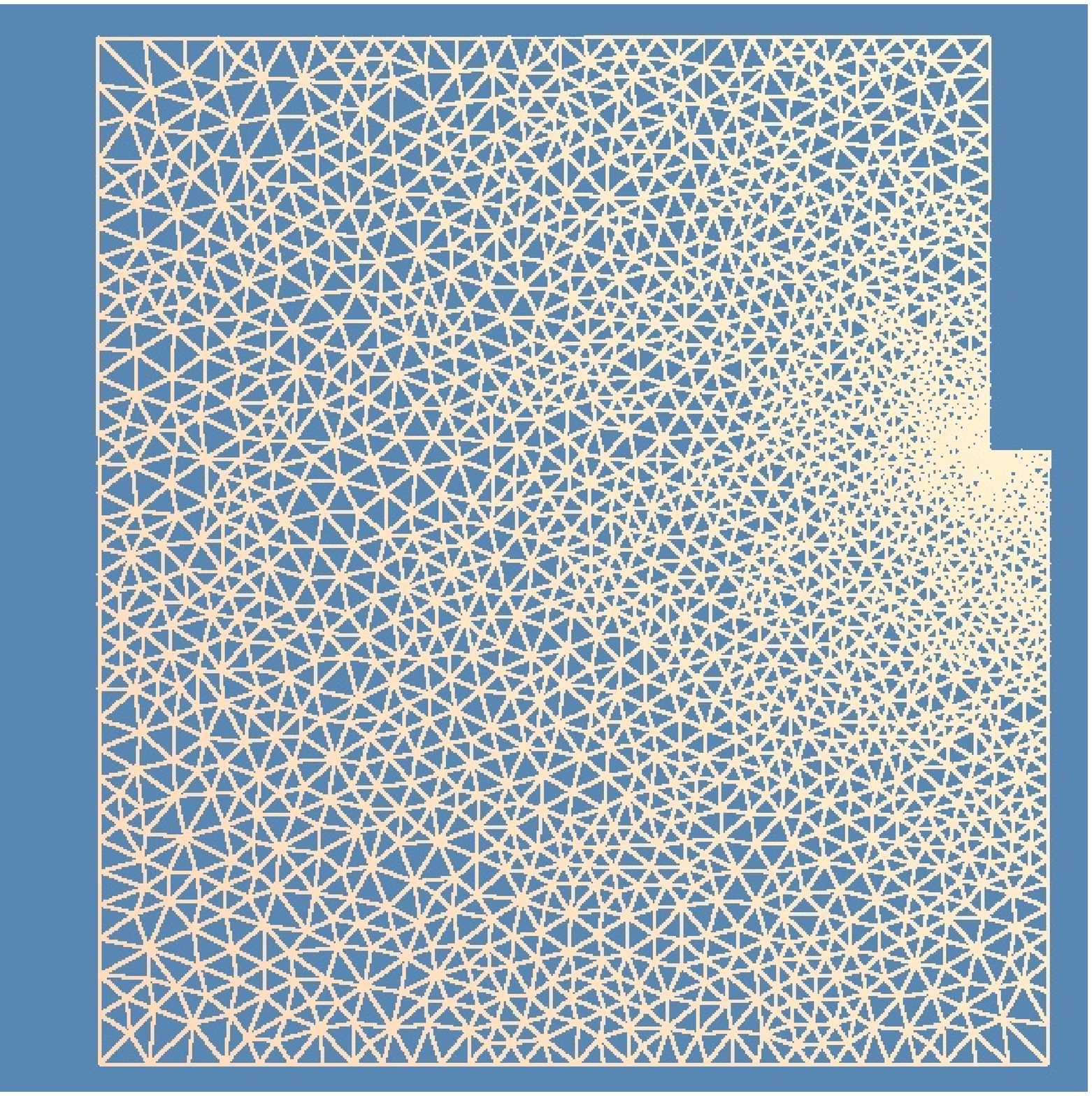}\\
(a)&(b)& (c)& (d)\\
\end{tabular}
\caption{Holomorphic quadratic differentials on a pentagon.  (a) and (b) show $[0.2(\phi_1')^2 + 0.8 (\phi_2')^2]dz^2$, (c) and (d) show $[-0.2(\phi_1')^2 + 1.2 (\phi_2')^2]dz^2$.\label{fig:pentagon_quad_diff}}
\end{figure*}

\begin{figure*}[h!]
\centering
\begin{tabular}{cccc}
\includegraphics[width=0.25\textwidth]{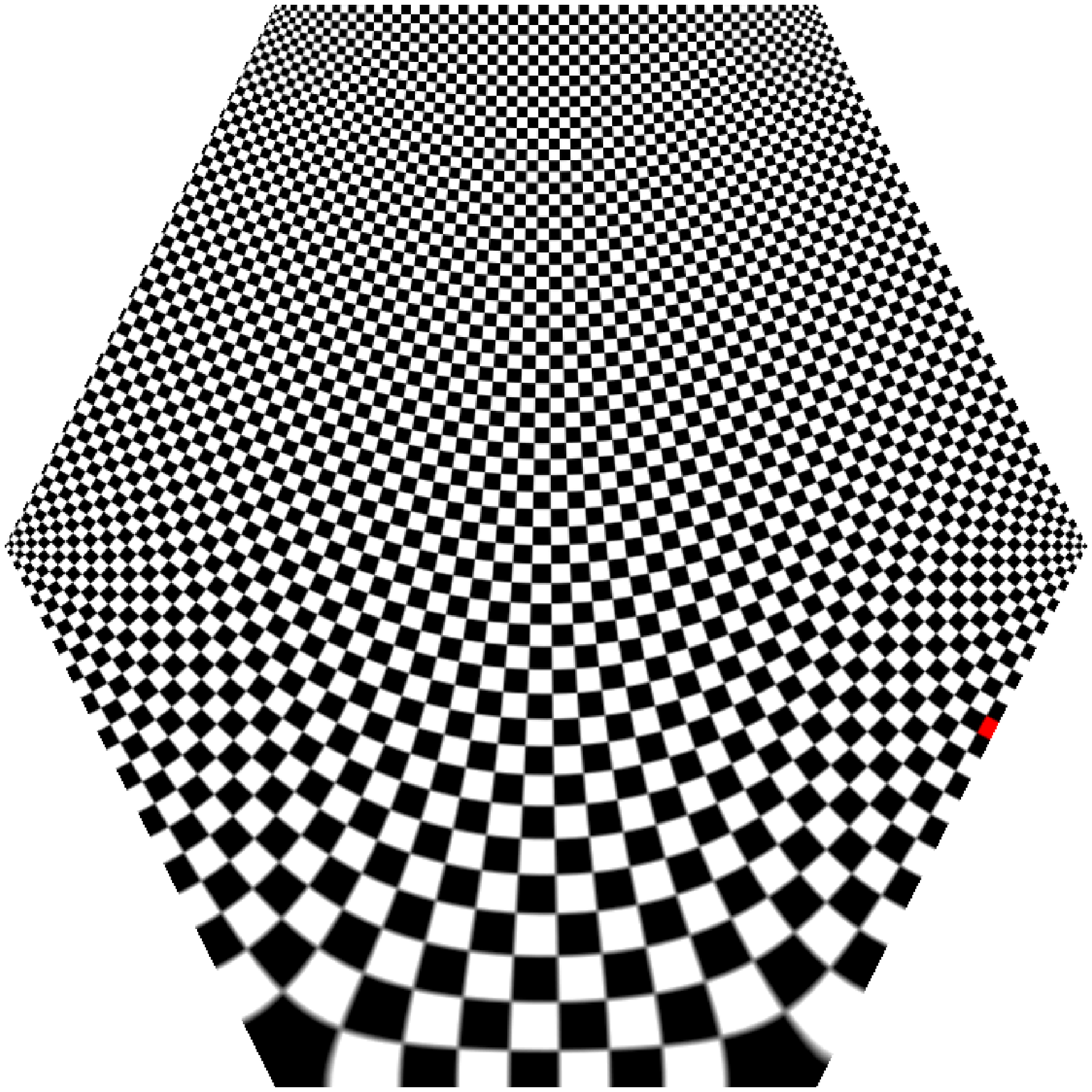}&
\includegraphics[width=0.25\textwidth]{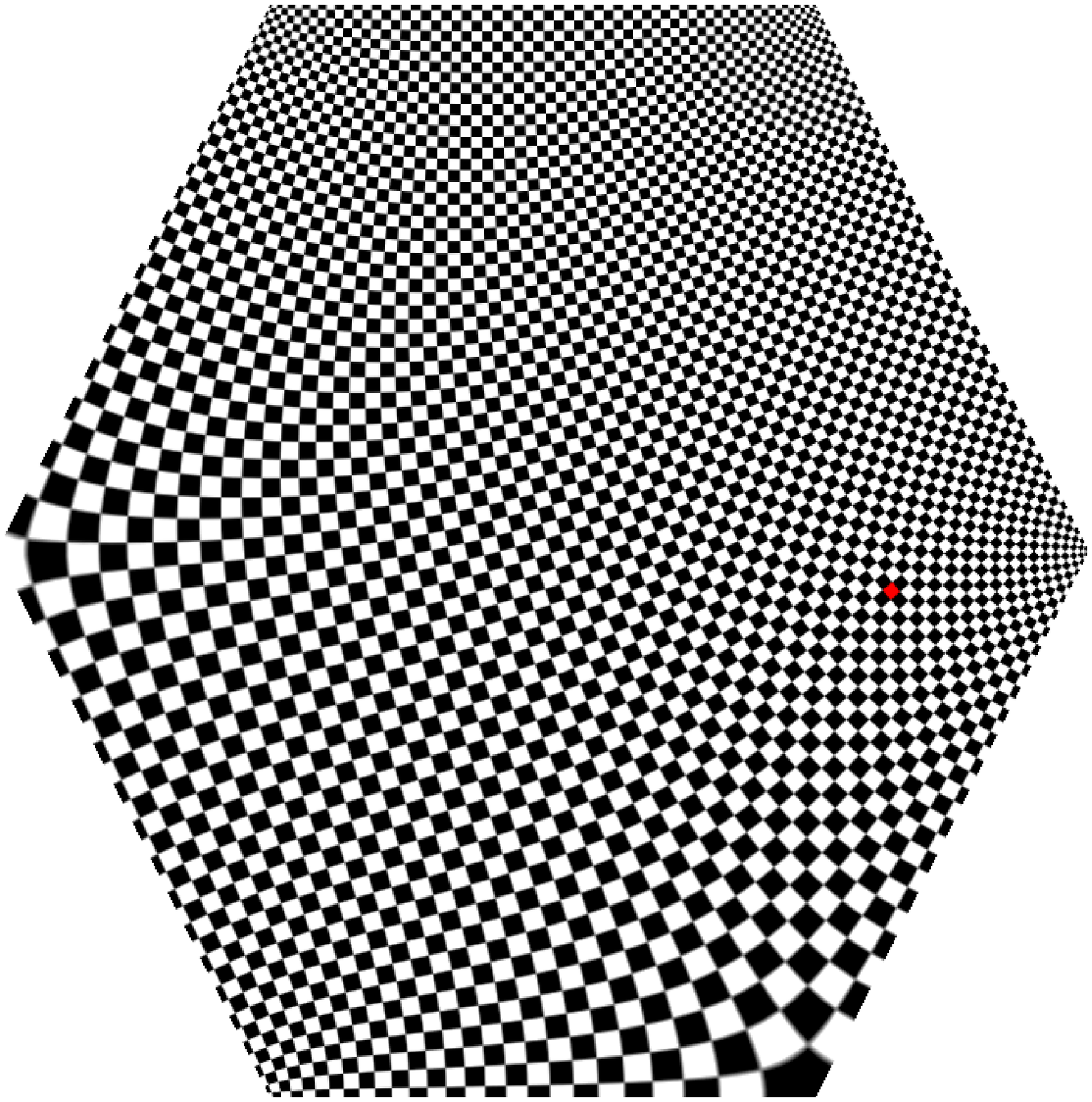}&
\includegraphics[width=0.25\textwidth]{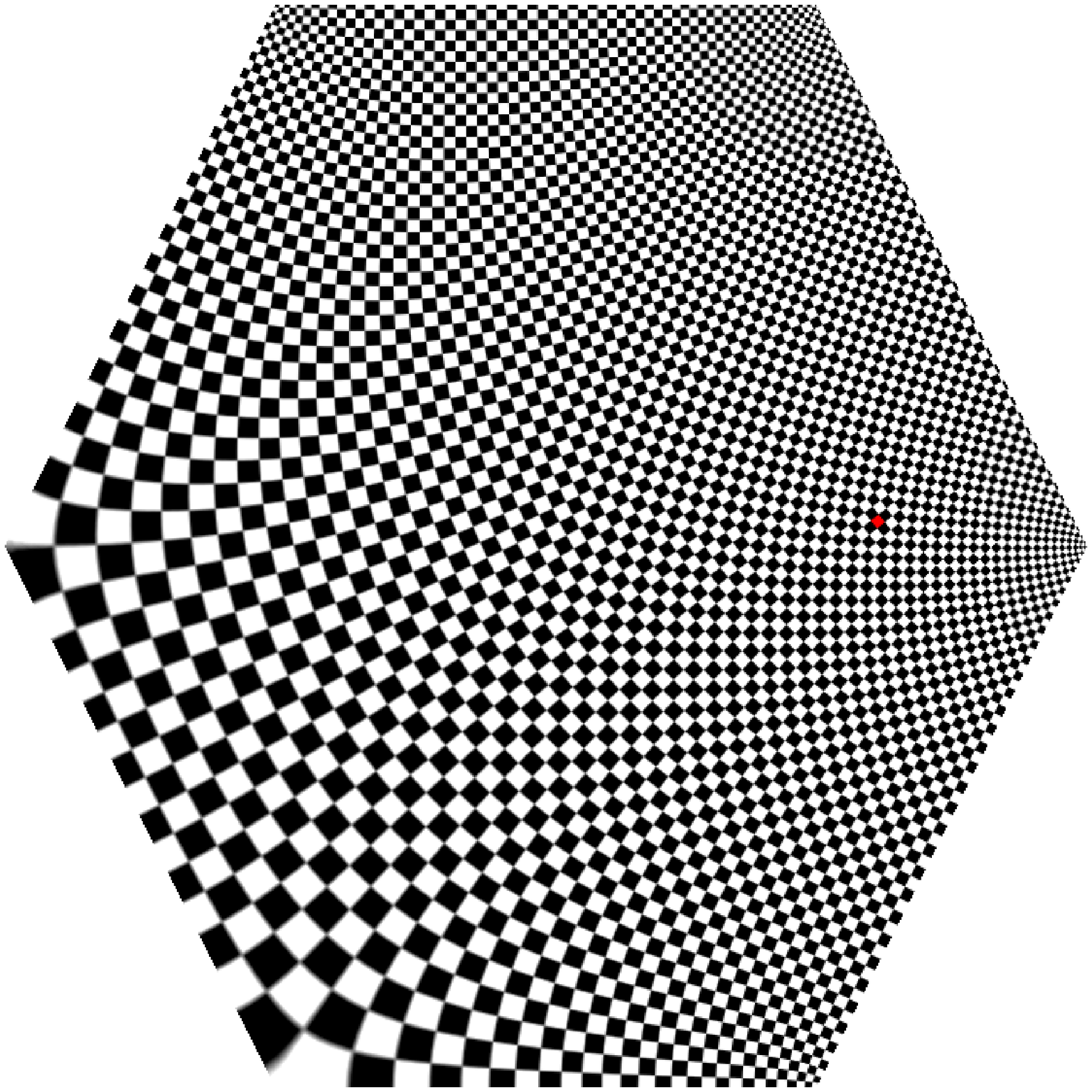}\\
\includegraphics[width=0.25\textwidth]{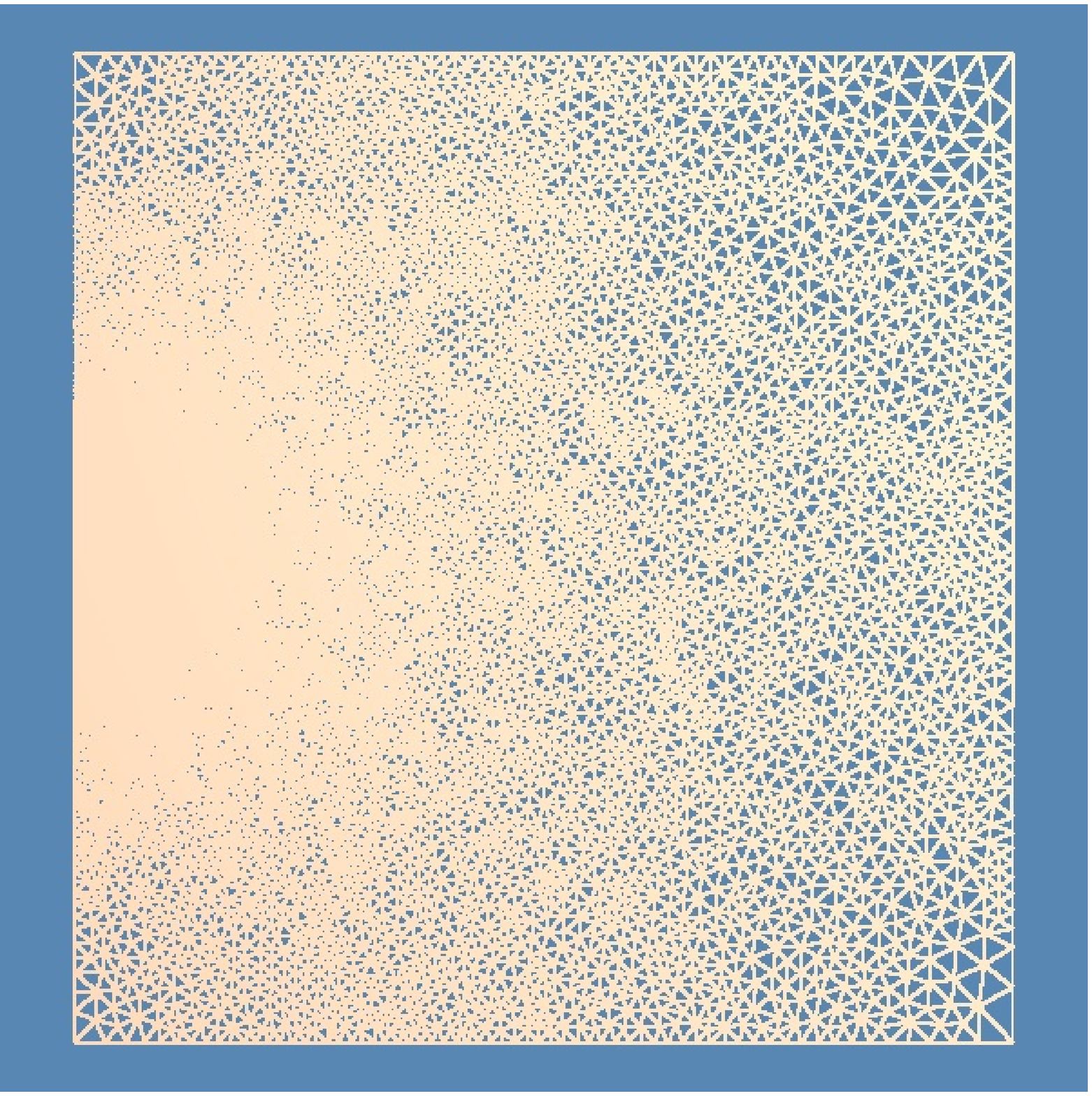}&
\includegraphics[width=0.25\textwidth]{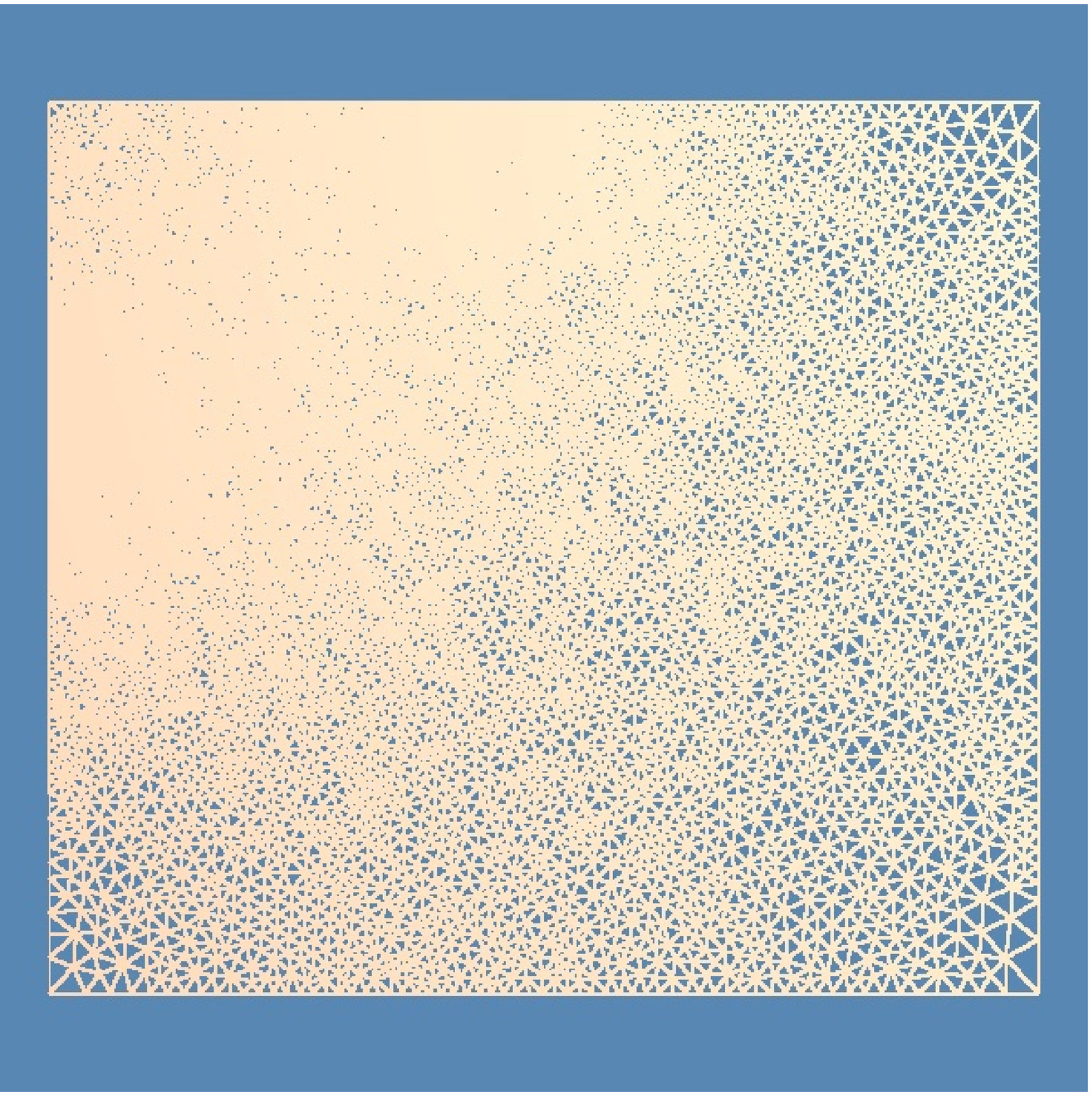}&
\includegraphics[width=0.25\textwidth]{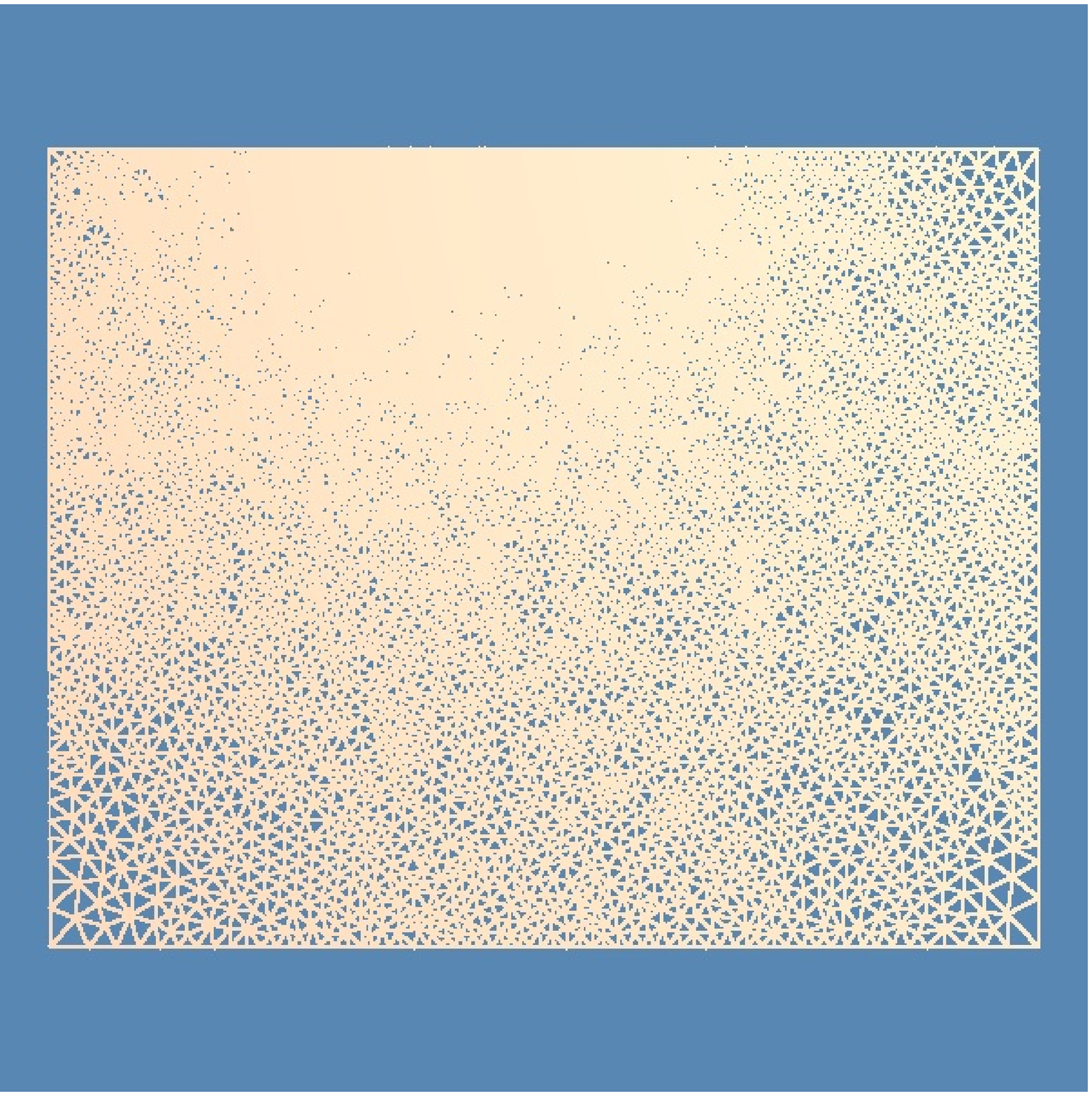}\\
\end{tabular}
\caption{Holomorphic quadratic differential bases on a hexagon. \label{fig:hexagon_quad_diff_base}}
\end{figure*}

\begin{figure*}[h!]
\centering
\begin{tabular}{cccc}
\includegraphics[width=0.25\textwidth]{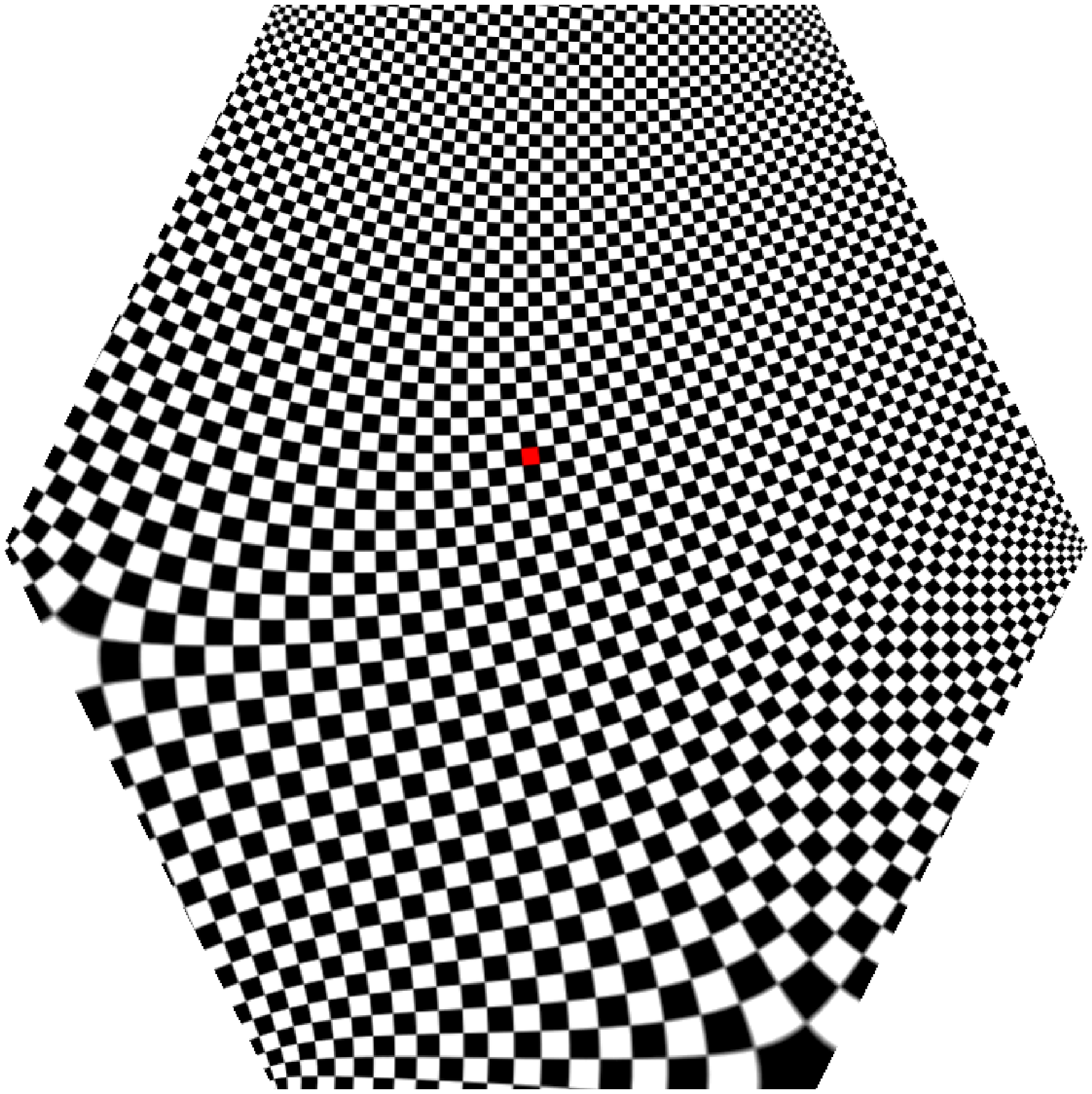}&
\includegraphics[width=0.25\textwidth]{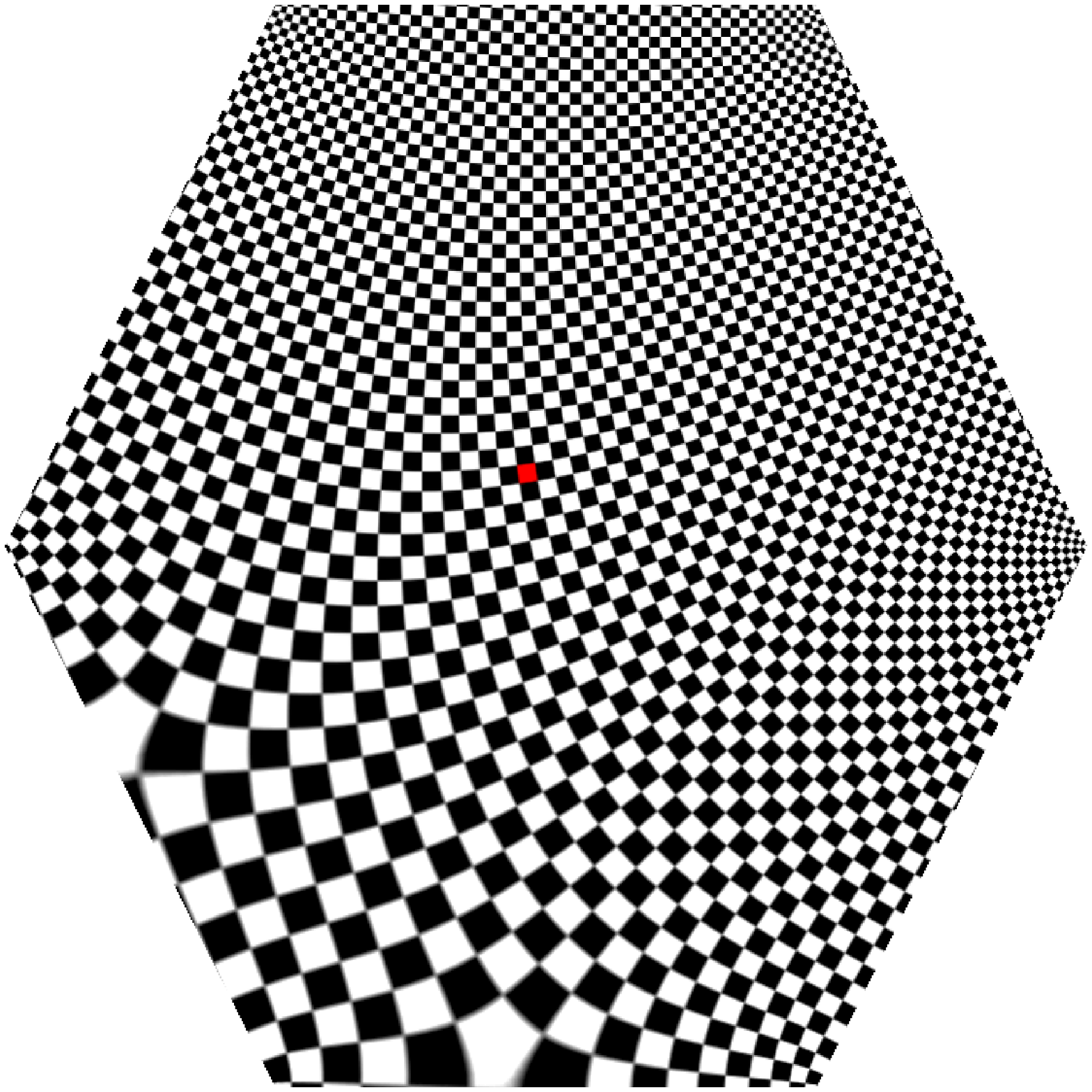}&
\includegraphics[width=0.25\textwidth]{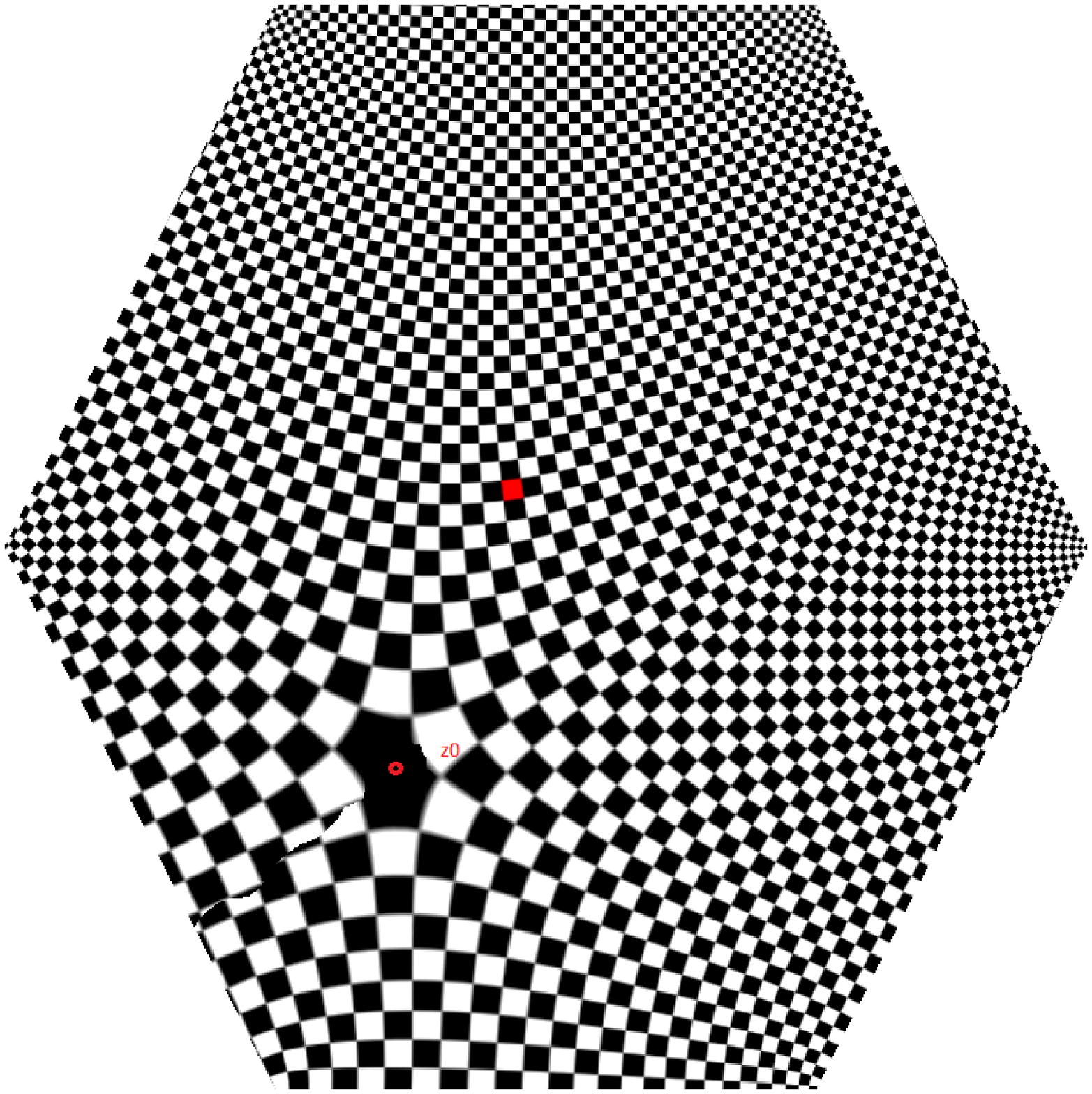}\\
\includegraphics[width=0.25\textwidth]{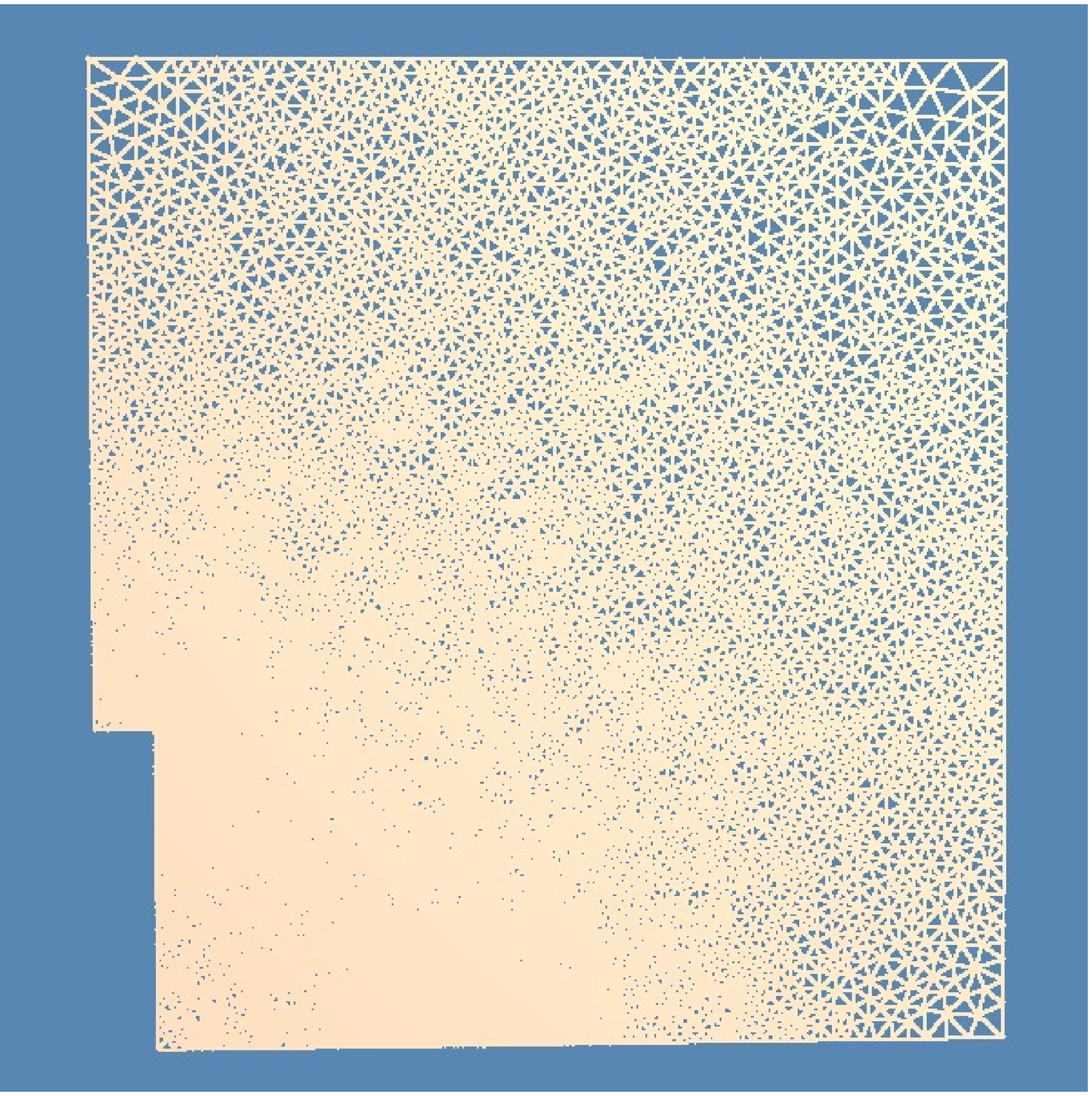}&
\includegraphics[width=0.25\textwidth]{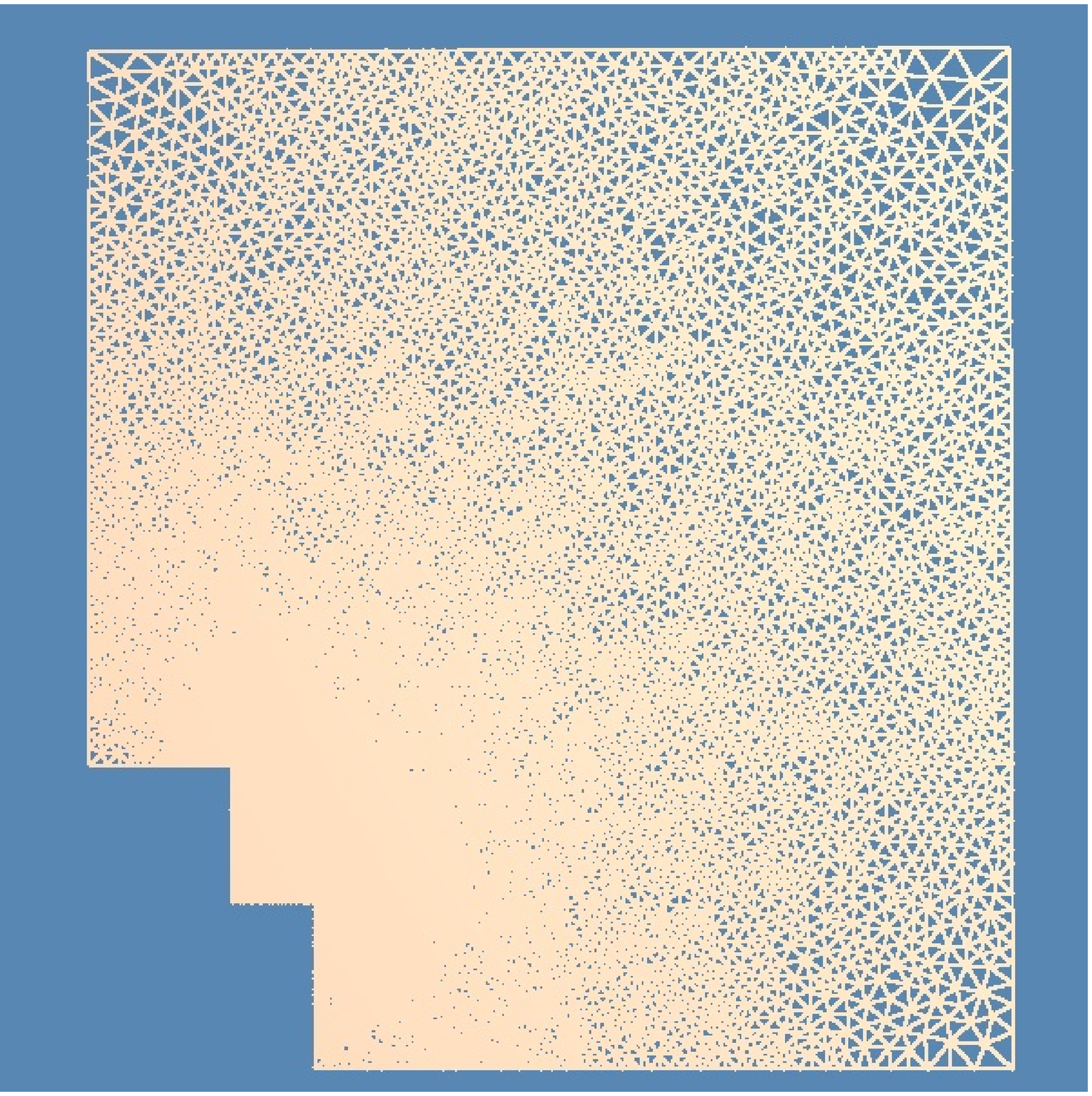}&
\includegraphics[width=0.25\textwidth]{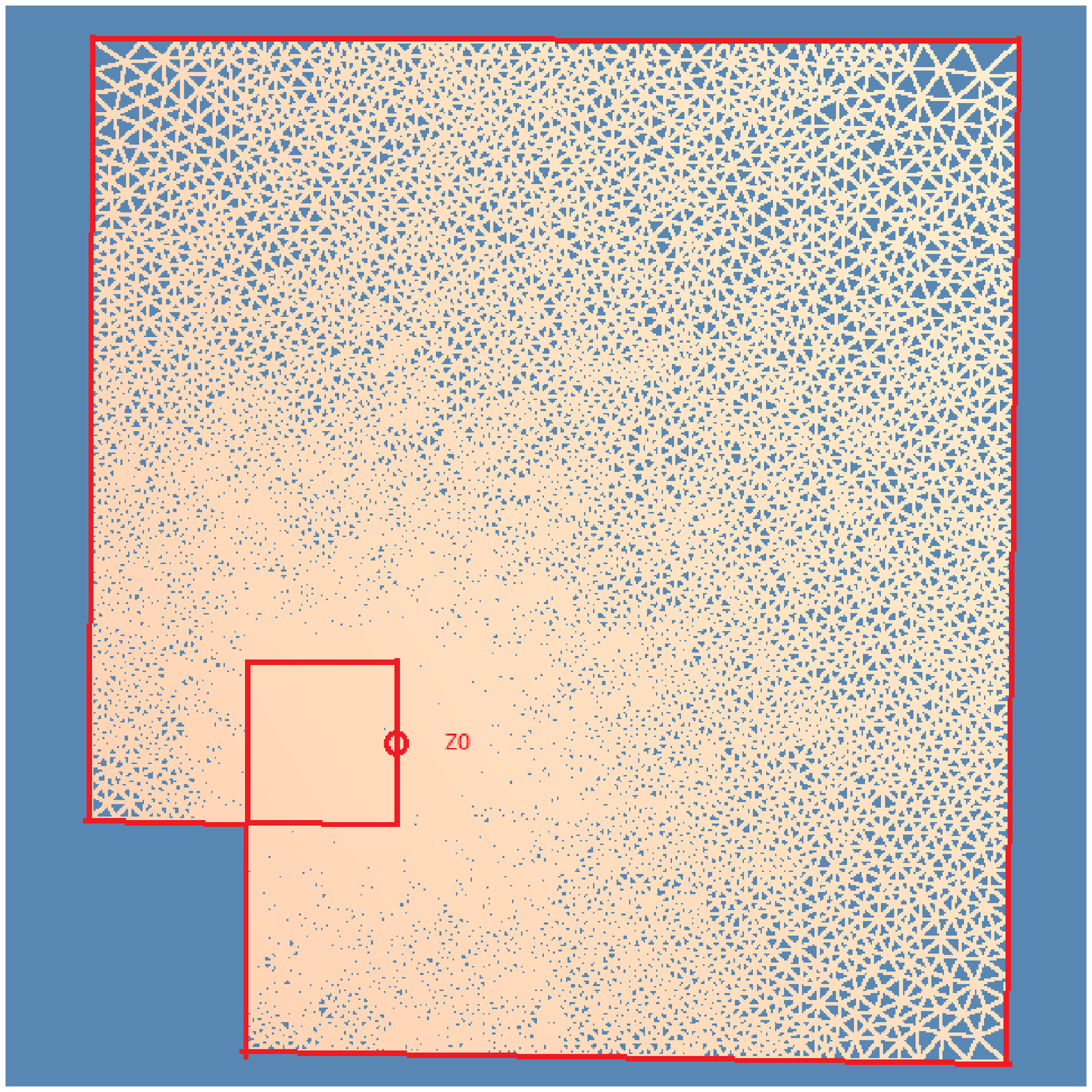}\\
\end{tabular}
\caption{Holomorphic quadratic differentials on a hexagon. \label{fig:hexagon_quad_diff}}
\end{figure*}

\begin{definition}[Holomorphic quadratic differential]
 A holomorphic quadratic differential on a Riemann surface $R$ is an assignment of a function $\phi_{i}(z_{i})$ on each chart $z_{i}$ such that if $z_j$ is another local coordinate, then $\phi_{i}(z_{i})=\phi_{j}(z_{j}){(\frac{dz_{j}}{dz_{i}})}^{2}$.
\end{definition}
We will denote the space of such differentials on $R$ as $A(R)$. By the Riemann-Roch theorem, the complex dimension of this  vector space for a genus $g$ closed compact surface with $n$ punctures is $3g-3+n$.

The Riemann surfaces of primary importance to us are the punctured Riemann sphere and the unit disk. For the unit disk, there is only one chart $z$, and therefore any function holomorphic in the interior of the disk can be viewed as a quadratic differential (the transition condition is vacuous). For $R = \hat{\mathbb{C}} \setminus \{0,1,\infty,z_{1},...z_{n-3}\}$ (the Riemann sphere with $n$ punctures),
\begin{equation}
\phi_{i}(z) = \frac{z_{k}(z_{k}-1)}{z(z-1)(z-z_{k})}, \ \ \ 1\leq k \leq n-3,
\end{equation}
form a basis of $(n-3)$ dimensional complex vector space $A(R)$.

Another vector space of importance to us is the space of polygon differentials. Let $P$ be a polygon in the plane, normalized so that $0$, $1$ and $\infty$ are three vertices of $P$. Suppose $\phi_k: P \to D_k$ is the conformal mapping, where $D_k$ is a planar rectangle, such that $\phi_k$ maps $\{0,1,\infty, z_k\}$ to the four corners of the rectangle $D_k$. Then
\[
   \{ (\phi_1')^2dz^2, (\phi_2')^2dz^2,\cdots, (\phi_{n-3}')^2dz^2 \}
\]
form the bases of $A(R)$. As shown in Figure~\ref{fig:pentagon_quad_diff_base}, the Riemann surface $R$ is a pentagon with vertices $\{z_1,z_2,z_3,z_4,z_5\}$, $\phi_1$ maps $R$ to planar rectangles $R_1$, such that $\{z_1,z_2,z_3,z_4\}$ are mapped to four corners. The checkerboard texture on $R_1$ is pulled back to $R$ and shown in (a). Similarly, $\phi_2$ maps $\{z_1,z_2,z_3,z_5\}$ to a rectangle $R_2$. Then $\{(\phi_1')^2dz^2, (\phi_2')^2dz^2\}$ form the bases of all holomorphic quadratic differentials on the pentagon. Figure~\ref{fig:pentagon_quad_diff} shows the linear combinations of these bases. Figures~\ref{fig:hexagon_quad_diff_base} and \ref{fig:hexagon_quad_diff} show the bases and certain linear combinations of the bases on a hexagon, respectively.

An excellent book for studying Quadratic differentials in further detail is \cite{strebel}.

\subsection{Classical theorems used in our construction}
\subseclabel{ingredients}

We start with a theorem which explains the dependence of a Beltrami coefficient to the solution of its Beltrami equation.

\paragraph*{Mapping Theorem}
[\cite{Gardiner}, Theorem $1$, Page $10$] Let $\mu(z)$ be a measurable complex-valued function defined on a domain $\Omega$ for which $||\mu||_{\infty} = k < 1$. Consider the Beltrami equation,

\begin{equation} \label{be}
f_{\bar{z}}(z)=\mu (z) f_{z}(z).
\end{equation}

\begin{theorem}
\thmlabel{mapping}
Equation~\ref{be} gives a one to one correspondence between the set of quasiconformal homeomorphisms of $\hat{\mathbb{C}}$ that fix the points $0,1$ and $\infty$ and the set of measurable complex-valued  functions $\mu$ on $\hat{C}$ for which $||\mu||_{\infty} <1$. Furthermore, the normalized solution $f^{\mu}$ of Equation~\ref{be} depends holomorphically on $\mu$ and for any $r>0$ there exists $\delta >0$ and $C(r) >0$ such that

\begin{equation} \label{stb}
|f^{t \mu}(z) - z - t V(z)| \leq C(r)t^{2} \ for \ |z|<r \ \text{and}\ |t|<\delta,
\end{equation}
where
\begin{equation} \label{derivative}
V(z) = - \frac{z(z-1)}{\pi}\int \int_{\mathbb{C}} \frac {\mu(\zeta)d \xi d \eta}{\zeta (\zeta-1)(\zeta - z)},
\end{equation}

and $\zeta = \xi + i \eta $.

\end{theorem}

\paragraph{Composition of Quasiconformal Maps}

Let $\mu$, $\sigma$ and $\tau$ be the Beltrami coefficients of quasiconformal maps $f^{\mu}$, $f^{\sigma}$ and $f^{\tau}$ with $f^{\tau} = f^{\sigma} \circ (f^{\mu})^{-1}$. Then

\begin{equation}\eqlabel{formula}
 \tau = \left( \frac{\sigma - \mu}{1 - \bar{\mu}\sigma} \frac{1}{\theta}\right) \circ (f^{\mu})^{-1},
 \end{equation}

where $p = \frac{\partial}{\partial z} f^{\mu}(z)$ and $\theta = \frac{\bar{p}}{p}$. In particular, if $f^{\sigma}$ is the identity, that is, if $\sigma = 0$, then
\[ \tau = - \left(\mu \frac{p}{\bar{p}}\right) \circ (f^{\mu})^{-1}. \]

The following lemma relates infinitesimally trivial Beltrami coeffcients to globally trivial ones.

\begin{lemma}\lemlabel{variationallemma}[Variational lemma][\cite{Gardiner}, Theorem $6$, Page $140$]
$\mu$ is an infinitesimally trivial Beltrami differential if, and only if, there exists a curve $\sigma_{t}$ of trivial Beltrami differentials for which $\sigma_{t}(z)=t \mu_{z}+O(t^{2})$ uniformly in $z$.
\end{lemma}

\paragraph{Teichm\"{u}ller contraction}\label{tcontract} The principle of Teichm\"{u}ller contraction states that given a Beltrami coefficient $\mu$, its distance to the globally extremal $\mu^{*}$ is of the same order as its distance to the infinitesimally extremal $\upsilon$. For a full statement and proof of the principle, see \cite{Gardiner}, Theorem $10$, page $103$.

We will restate the part of the principle relevant to us. Let $k_{0} = ||\mu^{*}||_{\infty}$ be the dilatation of the extremal Beltrami coefficient in the same global class as $\mu$, and let $\upsilon$ be the infinitesimally extremal Beltrami coefficient in the infinitesimal class of $\mu$. Fix $0 < k_{1} < 1$. Then

\begin{equation}\eqlabel{tcontract}
\frac{||\mu||_{\infty}-k_{0}}{4} \leq \frac{2}{(1-k_{1})^{2}}(||\mu||_{\infty}-||\upsilon||_{\infty}) \leq \frac{2}{(1-k_{1})^{4}}(||\mu||_{\infty}-k_{0}).
\end{equation}

\paragraph{Hamilton-Krushkal, Reich-Strebel condition}

\begin{theorem}\thmlabel{hkrs}[Hamilton-Krushkal, Reich-Strebel necessary-and-sufficient condition for extremality]
A quasiconformal map $f$ has minimal dilatation in its Teichm\"{u}ller class if and only if its Beltrami coefficient $\mu$ is extremal in its infinitesimal class.
\end{theorem}

\section{Appendix for the discretization of the procedure}
\seclabel{discrete_appendix}

\subsection{Formula for the integral of $\phi_{i}$ on a triangle $t_{j}$}
\subseclabel{formula_integral}

Let D be the triangle whose vertices are $\al$, $\be$, $\ga$ (in that order). Then
\begin{gather}
\int _{D} \frac{1}{z(z-1)(z-a)} = \int_{D} \left( \frac{1}{az} + \frac{1}{(1-a)(z-1)} + \frac{1}{a(a-1)(z-a)} \right ) \nonumber \\
 = \int _{0} ^{1} \int _{0}^u (\frac{1}{a(\al+u(\be-\al)+v(\ga - \be))} + \frac{1}{(1-a)(\al-1+u(\be-\al)+v(\ga - \be))} + \nonumber \\ \frac{1}{a(a-1)(\al-a+u(\be-\al)+v(\ga - \be))} dv du \nonumber \\
 = \frac{1}{\ga-\be}\displaystyle \left[ \frac{1}{a}I_1 + \frac{1}{1-a}I_2 + \frac{1}{a(a-1)}I_3 \right] \nonumber\\
 I_1 = \frac{1}{\ga-\al}\left( \ga \ln \ga - \ga - \al \ln \al + \al \right)-\frac{1}{\be-\al}\left( \be \ln \be - \be - \al \ln \al + \al \right )\nonumber \\
 I_2 = \frac{1}{\ga-\al}((\ga-1) \ln (\ga-1) - \ga - (\al-1) \ln (\al-1) + \al)\nonumber \\ -\frac{1}{\be-\al}((\be-1) \ln (\be-1) - \be - (\al-1) \ln (\al-1) + \al))\nonumber \\
 I_3 = \frac{1}{\ga-\al}((\ga-a) \ln (\ga-a) - \ga - (\al-a) \ln (\al-a) + \al)\nonumber \\ -\frac{1}{\be-\al}((\be-a) \ln (\be-a) - \be - (\al-a) \ln (\al-a) + \al)) \nonumber \\
\end{gather}

\subsection{Details of the algorithm INEXT}
\subseclabel{inextdetails}
$\mathcal{P}(\mu)$ is the solution to the following program.
\begin{program}\label{continuousprogram}

\begin{center}
\[\min \ \ \ ||\nu||_{\infty}\]
\[\text{subject to}: \ \int_{R} \nu \phi_{i} = \int_{R} \mu \phi_{i} \ \forall i \in \{1,2...n-3\}\]
\end{center}
\end{program}

Program~\ref{continuousprogram} is an $L_{\infty}$ norm minimization subject to certain constraints. We will solve the above program when $\nu$ ranges over all piecewise constant Beltrami differentials. Let $\{\nu_{i}\}_{i=1}^{T}$, where $T$ is the number of triangles in the triangulation, be a basis of piecewise constant Beltrami differentials. 

Consider the $j$th triangle $t_j$; the integral of any basis element $\phi$ in Equation~\eqreff{puncturediff} can be computed analytically, and in a preprocessing step, we compute the matrix $A$ where $((a_{ij})) = \int_{t_{j}} \phi_{i}$. 

We write $\nu = \sum \lambda_{i} \nu_{i}$, and represent $\nu$ as a vector $\lambda = (\lambda_1,\lambda_2,\cdots, \lambda_T)$. Each $\lambda_k$ is a complex number; separating into real and imaginary parts we get  $\lambda_k = \lambda_{rk} + i\lambda_{ik}$. Let the analogous vector for $\mu$ be $\lambda^{'}$. Then the above constraint can be written as $A \lambda = b$, where $b = A \lambda^{'}$. Using the above separation into real and imaginary parts for the matrix $A = A_{r} + i A_{i}$ and $b = b_{r} +ib_{i}$, this is equivalent to

\begin{definition}[Constraints]
\begin{equation} \label{cons1}
 A_{r}\lambda_{r} - A_{i}\lambda_{i} = b_{r}
\end{equation}
\begin{equation} \label{cons2}
A_{i}\lambda_{r} + A_{r}\lambda_{i} = b_{i}
\end{equation}
\end{definition}

We introduce another variable $z \in \mathbb{R}^{+}$, making the number of variables $(2T+1)$. Let the vector of unknowns be $\beta = (\lambda_{r1},\cdots,\lambda_{rT},\lambda_{i1},\cdots,\lambda_{iT},z)$.

\begin{program}\label{discreteprogram}
\begin{center}
\[ \min \ \ \ z \]
\begin{eqnarray}
\text{subject to }&:& \text{Constraints }\ref{cons1} \text{ and }\ref{cons2},\notag \\
& \text{and }& \lambda_{rj}^{2} + \lambda_{ij}^{2} - z \leq 0 \  \ \ \forall 1 \leq j \leq T \notag
\end{eqnarray}
\end{center}
\end{program}

The last constraint uses the fact that the solution $\nu^{*}$ to Program~\ref{continuousprogram} is of Teichm\"{u}ller form. The objective function of Program~\ref{discreteprogram} is linear in the unknown variables. Constraints~\ref{cons1}~and~\ref{cons2} are also linear. The last set of constraints can be written as $\beta^{t} P_{j} \beta - z  \leq 0$, where $P_{j}$ is a $(2T+1)$ matrix of zeroes with its $(j,j)^{th}$ and $(T+j,T+j)^{th}$ entry being $1$. $P_{j}$ has all but two eigenvalues as $0$, and two eigenvalues are $1$, implying that it is positive semi-definite.

Although solving a quadratically constrained quadratic program in general is NP-Hard, positive semi-definite instances of it are polynomial time solvable. Numerical solvers for these programs have been vastly studied, and efficient implementations exist. We refer the reader to Page~$42$ of \cite{pardalos2002handbook} for a complete reference. 

This completes the proof for \lemref{def_inext}.

\subsection{Holomorphic quadratic differentials on polygons}
\subseclabel{qdonpoly}
Suppose $P$ is a polygon with vertices $\{z_0,z_1,\cdots, z_{n-1}\}$. For each $2<k<n$, there exists a unique conformal map $\phi_k$, which maps $P$ to a rectangle $\R = [0,1]\times[0,h]$, and maps $z_0,z_1,z_2, z_k$ to the four vertices of $\R$. Then $\{(\phi_k')^2dz^2\}$ form the basis of holomorphic quadratic differentials on $P$.

All the above proofs can be modified analogously for the polygon mapping problem with minimal effort.

\end{document}